\documentclass[12pt]{amsart}
\usepackage{graphicx,calc,amssymb}

\newtheorem{itheorem}{Theorem}
\newtheorem{theorem}{Theorem}[subsection]
\newtheorem{corollary}[theorem]{Corollary}
\newtheorem{proposition}[theorem]{Proposition}
\newtheorem*{remarks}{Remarks}
\newtheorem{remark}[theorem]{Remark}
\newtheorem{rmk}[theorem]{}
\newtheorem{lemma}[theorem]{Lemma}
\newtheorem{defi}[theorem]{Definition} 

\newcommand{\X}{{}{\ensuremath{X_{ij}}}}

\newcommand{\Y}{{}{\ensuremath{Y_{ij}}}}
\newcommand{\C}{\operatorname{\mathcal{C}}}
\newcommand{\f}{{}{\ensuremath{f_{ij}}}}

\newcommand{\GLn}{\operatorname{GL_n}}
\newcommand{\SL}{\operatorname{SL_n}}
\newcommand{\GL}{{}{\ensuremath{\GLn(C)}}}

\newcommand{\bx}{{}{\ensuremath{\mathbf{X}}}}
\newcommand{\ba}{{}{\ensuremath{\mathbf{\alpha}}}}
\newcommand{\bb}{{}{\ensuremath{\mathbf{\beta}}}}
\newcommand{\Lie}{\operatorname{Lie}}
\begin{document}

\title{On a generic inverse differential Galois problem for GL$_\mathbf n$}

\author{Lourdes Juan}
\address{Mathematical Sciences Research Institute\\
1000 Centennial Dr.\\
Berkeley, California 94720\\
USA}
\email{ljuan@msri.org}
\urladdr{www.msri.org/people/members/ljuan/}

\subjclass{Primary 12H05; Secondary 12F12, 12Y05, 20G15}

\date{\today}
\thanks{ This paper was written while the 
 author held an NSF-funded  postoctoral fellowship at the MSRI. The author wishes to
 thank the MSRI for its hospitality.}
\keywords{Constants, differential field, differential Galois group,
 differential specialization, inverse differential Galois problem, generic, multivariable division algorithm,
 new constant,  Picard-Vessiot  extension, wronskian. }

\begin{abstract} In this paper we construct a generic Picard-Vessiot extension  for the general linear groups.
In the case when the   differential base field has finite transcendence degree over its
field of constants 
we provide necessary and sufficient conditions for  solving the inverse differential Galois problem for this groups
via specialization  from  our generic extension.
\end{abstract}

\maketitle
\markboth{LOURDES JUAN}{ON A GENERIC INVERSE DIFFERENTIAL GALOIS PROBLEM FOR  GL$_\mathbf n$}

\section*{Introduction}
\label{sec:intro}
 
 Le  $F$ be a differential field of characteristic zero 
 with algebraically closed field of constants $C$. In this paper we
  give an affirmative answer, for the group $\GL$, to the following 

{\bf Generic Inverse Differential Galois Problem}: \emph{For a connected algebraic group $G$ over $C$
find a generic Picard-Vessiot  extension of $F$ with differential Galois group $G$.}

By \emph{generic extension} we mean a Picard-Vessiot  extension of a \emph{generic field} that contains $F$ and
such that every Picard-Vessiot  extension of $F$ for $G$ in the usual sense can be obtained from the
generic one by specialization.

We point out that any such specialization will provide a solution to the inverse differential
Galois problem in the usual sense, namely, to
 determine,  given $F$ and $C$ as above, and a linear
algebraic group $G$ over $C$,  what differential
field extensions $E\supset F$ are Picard-Vessiot  extensions   with differential
Galois group $G$ and, in particular, whether there are any. 

We use the   terminology   of A. Magid's  
book \cite{M}. In \cite{M} the reader  may also find definitions
and  proofs of some results from differential Galois theory that
will be recalled. 

Our construction of
a generic Picard-Vessiot  extension of $F$ with group $\GL$ may be summarized as follows: Let $F\{\Y\}$ 
be the ring of differential polynomials over $F$ in the differential indeterminates
$\Y$, $1\le i, j\le n$. Our   \begin{emph}{generic base field}\end{emph} 
  will be the differential quotient field 
$F\langle \Y\rangle$  of $F\{\Y\}$.  Let $\X$, $1\leq i,j\leq n$,
be algebraically independent over $F\langle \Y\rangle$. Extend the derivation on $F\langle \Y\rangle$ 
to the polynomial
ring $F\langle \Y\rangle[\X]$ by letting $D(\X)=\sum_{\ell=1}^n Y_{i\ell}X_{\ell j}$.
This derivation extends in a natural way to the quotient field $F\langle \Y\rangle(\X)$
 (note that 
$F\langle \Y\rangle(\X)$ is  also the
function field of the group obtained from $\GL$ 
by extending  scalars  from $C$ to $F\langle\Y\rangle$). 
Then 
\begin{itheorem}\label{ige} The differential field extension
 $F\langle \Y\rangle(\X)\supset F\langle \Y\rangle$ is a generic 
 Picard-Vessiot extension of $F$ with differential Galois group $\GL$.\end{itheorem}

The main step in proving Theorem ~\ref{ige} is to show that 
$F\langle \Y\rangle(\X)\supset F\langle \Y\rangle$ is a no-new-constant extension. A
new constant in $F\langle \Y\rangle(\X)$ must be the quotient of two relatively 
prime Darboux polynomials, that is, two polynomials $p_1,p_2\in F\{\Y\}[\X]$ 
that satisfy $D(p_i)=qp_i$ 
for some $q\in F\{\Y\}[\X]$. We use Gr\"obner basis techniques to show that the only Darboux
polynomials in  $F\{\Y\}[\X]$ are those of the form $\ell\det[\X]^a$, where $\ell\in F$ and
$a\in \mathbb N$. This implies that there are no two such relatively prime Darboux polynomials
and, therefore, no new constants. 

To show that the extension is generic we let $E\supseteq F$ be any Picard-Vessiot  extension of $F$ with
group $\GL$. Then $E$ is isomorphic to $F(\X)$, 
the function field of the group obtained from $\GL$ by extending scalars from
$C$ to $F$,
as a $\GL$-module and as an $F$-module, and there are $\f\in F$
such  that the derivation on $E$ is given by $D_E=\sum_{\ell=1}^n f_{i\ell}X_{\ell j}$. 
In this situation, the Picard-Vessiot  extension
$F(\X)\supseteq F$ is obtained from the Picard-Vessiot  extension 
$C\langle \f\rangle(\X)\supseteq C\langle \f\rangle$ by 
extension of scalars from $C$ to $F$. Therefore, any Picard-Vessiot extension
$E\supset F$ can be obtained from $F\langle \Y\rangle(\X)\supset F\langle \Y\rangle$
via the specialization $\Y\mapsto \f$. That is, $F\langle \Y\rangle(\X)\supset F\langle \Y\rangle$
is a generic Picard-Vessiot extension of $F$ for $\GL$.

Now, suppose that $F$ has finite transcendence degree over $C$ say, $$F=C(t_1,\dots,t_m)[z_1,\dots,z_k],$$
where the $t_i$ are algebraically independent over  $C$ and the $z_i$ are algebraic over $C(t_1,\dots,t_m)$. 
Consider the differential field  $F(\X)$ with derivation  $D_E$ as
above. Let $\C$ denote its field
 of constants.  Let  $F\{\Y\}[\X]$ be the differential ring  with the derivation  defined  above.
For $k\geq 1$ let $\mathbb T_k$ denote the set of monomials in the both the $t_i$ and $\X$
 of total degree less than or equal to $k$.
Fix a term order on the set $\mathbb T$ of  monomials in the $t_i$ and the
 $\X$ and let $W_k(\Y)$ denote the wronskian 
 of $\mathbb T_k$
relative to that order (the order will only affect the wronskian by a sign). The following
theorem summarizes our specialization results:

\begin{itheorem}\label{isp}  
 $F(\X)\supset F$ is a Picard-Vessiot  extension
for $\GL$
 if and only if
all the wronskians $W_k(\Y)$ map to nonzero elements in $F(\X)$  via the
specialization $\Y\mapsto\f\in F$.
\end{itheorem}

The above condition on the wronskians means that
all the sets $\mathbb T_k$, for $k\geq 1$, are linearly independent over $\C$.
This is in turn  equivalent to the fact that the set of all the $t_i$ and
all the $\X$ are algebraically independent over $\C$.
Unfortunately, Theorem ~\ref{isp} gives infinitely many conditions. We do not know at present
how to use these conditions to effectively construct solutions to the inverse problem, and
this constitutes an interesting open problem. 

A specialization as in Theorem~\ref{isp},
however, is
known to exist by a result  of 
C. Mitschi and M. Singer  \cite{M-S2}. They give a constructive algebraic solution 
to the inverse problem for
all connected linear algebraic groups (and, in particular, for $\GL$) when
$F$ has finite transcendence degree over $C$. An interesting direction of research
in connection with the previous open problem
is to try to  fully describe all  possible solutions (isomorphic and
non-isomorphic) that may arise in this situation.

The work of Mitschi and Singer in \cite{M-S2} makes use of the logarithmic derivative and
an inductive technique developed by Kovacic  \cite{K1}, \cite{K2}, to lift a
solution to the inverse problem from $G/R_u$, where $R_u$ is the unipotent radical
of $G$, to the full group $G$. Using this machinery Kovacic proved that it is enough 
to find a solution to the inverse problem for reductive groups (observe that 
$G/R_u$ is reductive). In \cite{P}, van der Put explains and partly  proves  the 
 results in \cite{M-S2}.
 
In the introduction of  \cite{M-S2} the authors briefly review previous work on the inverse
problem such as 
results of Bialynicki-Birula in \cite{BB}, 
Kovacic \cite{K1}, \cite{K2}, Ramis \cite{ramis1}, \cite{ramis2}, Singer \cite{S2},  
Tretkoff and Tretkoff \cite{TT},
 Beukers and Heckman \cite{B-H},
Katz \cite{K}, Duval and Mitschi \cite{dm}, Mitschi \cite{M1}, \cite{M2},
Duval and Loday-Richaud \cite{dlr}, Ulmer and Weil \cite{uw}  and Singer and
Ulmer \cite{su}. A  more extensive survey on the inverse problem
can be found in M. Singer's \cite{S}. 

The constructive algebraic solutions to the inverse differential Galois
problem for connected algebraic groups
that are currently available are based on  Kolchin's Main Structure
Theorem for Picard-Vessiot  extensions (see Theorem ~\ref{sthm} below). In particular,
a corollary to this theorem (see Theorem ~\ref{stgln}) establishes that if 
$E\supset F$ is Picard-Vessiot and $G$ is, for example,  unipotent or
 solvable
or  $G=\GLn$ or  $G=\SL$, then
E is isomorphic as an $F$-module and as a $G$-module to the 
function field of the group $G_F$ obtained from $G$ by extension of
scalars from $C$ to $F$. Therefore, 
 to get a Picard-Vessiot  extension 
$E\supset F$ with group $G$ (if it exists) one can begin by taking 
$E$ to be the  function field of $G_F$ and then the problem
reduces to extending
the derivation from $F$ to $E$ in such a way that $E\supset F$ is Picard-Vessiot  
for that derivation. In this paper we use this approach for our construction.

The idea of tackling  the inverse problem by constructing  generic 
extensions is inpired by the works of E. Noether \cite{EN} for the Galois theory
of algebraic equations. Following her approach, L. Goldman  in \cite{goldman} introduced
the notion of a \begin{emph}{generic differential equation with group $G$}\end{emph}. Goldman explicitly
constructed a generic equation  with group $G$ for some groups. However, after specializing Goldman's
 equation  the  group of  the new equation obtained
is a subgroup of the original group. In order to solve the
inverse problem by this means, we need to keep the original group as
the group of the equation after specialization. 
Goldman's generic equation for $\GLn$ is equivalent to Magid's 
\begin{emph}{general equation of order $n$}\end{emph} (Example 5.26 in
\cite{M}).

More work in the spirit of Goldman's generic equation
came some years later in J. Miller's 
 dissertation  \cite{miller}. He  defined the notion of hilbertian differential
field  and gave a sufficient condition for the generic equation with 
group $G$ to specialize to an equation over such a field with group
$G$ as well. However, as pointed out
by Mitschi and Singer  in \cite{M-S2}, his condition was stronger 
than the analogous one for algebraic equations and this made the theory
especially difficult to apply for those groups that were not
already known to be Galois groups.

This paper contains the results of the author's Ph.D. dissertation \cite{J}. 
I wish to thank  my Ph.D. advisor  Andy Magid for the
many valuable research meetings that we had. I am also grateful to Michael Singer for 
 many  enlightening conversations on the inverse problem.

\section{Preliminaries}    

\subsection{Notation and some basic results from Differential Galois Theory}
We fix a differential field $F$ 
with algebraically closed field of
constants $C$. If $E\supseteq F$ is a differential
field extension then the group of differential automorphisms
of $E$ over $F$ is denoted by  $G(E/F)$.
 
 If $G$ is a linear algebraic group over $C$ and $K$ is an overfield of $C$
 we denote by $G_K$ the
 group obtained from $G$ by extending scalars from $C$ to $K$.
 
 We will show that the differential field $F\langle\Y\rangle(\X)$ to be defined
 in ~\ref{diffields} is a Picard-Vessiot extension  of $F\langle\Y\rangle$ with differential
 Galois group $\GL$. Note that \linebreak $F\langle\Y\rangle(\X)$  is the function
 field of $G_K$ with $G=\GL$ and $K=F \langle\Y\rangle$. The
 following two results provide the rationale for choosing 
 such function field for the generic Picard-Vessiot extension. The proofs
can be found in \cite{M}
  (Theorem 5.12 and Corollary 5.29 respectively). 
 
 \begin{theorem}[Kolchin Structure Theorem]\label{sthm}
 Let $E \supseteq F$ be a Picard-Vessiot extension, let $G\leq G(E/F)$
 be a Zariski closed subgroup  and let $T$ be the set of all $f$ in
 $E$ that satisfy a linear homogeneous differential equation over 
 $K=E^G$. Then $T$ is a finitely generated $G$-stable differential 
 $K$-algebra with function field $E$, and if $\overline K$ denotes the
 algebraic closure of $K$, then there is a $G$-algebra isomorphism
 $$\overline K\otimes_K T\rightarrow \overline K\otimes_C C[G].$$
 Note that $C[G]$ denotes the affine coordinate ring of $G$ and that the target of the
 above isomorphism is the affine coordinate ring of the group $G_{\overline K}$ obtained 
 from $G$ by extension of scalars from $C$ to $\overline K$.
 \end{theorem}
 
 \begin{theorem}\label{stgln}
 Let $E \supseteq F$ be a Picard-Vessiot extension, let $G\leq G(E/F)$
 be a Zariski closed subgroup  with $E^G=F$. Let $\overline F$ be an algebraic 
 closure of $F$, and suppose the Galois cohomology $H^1(\overline F/F,G(\overline F))$
 is a singleton. Let $T(E/F)$ be the set of all $f$ in $E$ that satisfy 
 a linear homogeneous differential equation over $F$.
 Then there are $F$- and $G$-isomorphisms $T(E/F)\rightarrow F[G_F]$ and
 $E\rightarrow F(G_F)$. In particular, this holds if $G$ is unipotent or solvable,
 or if $G=\GL$ or if $G=\SL$.\end{theorem}

 We will use the following characterization of Picard-Vessiot extension given
 in \cite{M} (Proposition 3.9):
 
 \begin{theorem}\label{cpve} Let $E \supseteq F$ be a differential field extension.
 Then $E$ is a Picard-Vessiot extension if and only if:
 \begin{enumerate}
 \item $E=F\langle V\rangle$, where $V\subset E$ is a finite-dimensional vector
 space over $C$;
 \item  There is a group $G$ of differential automorphisms of $E$ with $G(V)\supseteq V$
 and $E^G=F$;
 \item $E\supset F$ has no new constants.
 \end{enumerate}
 In particular, if the above conditions hold and if $\{y_1,\dots,y_n\}$ is a $C$-basis of
 $V$, then $E$ is a Picard-Vessiot extension of $F$ for the linear homogeneous
 differential operator
 $$L(Y)=\frac{w(Y,y_1,\dots,y_n)}{w(y_1,\dots,y_n)}$$
 where $w(-)$ denotes the wronskian determinant
 and $L^{-1}(0)=V$.\end{theorem}
 
 In our case the base field is $F\langle \Y\rangle$ and $G=\GL$. We want to show
 that  $F\langle \Y\rangle (\X)\supset F\langle \Y\rangle$ is Picard-Vessiot. We
 first show (Corollary ~\ref{constfy}) that the field of constants of $F\langle\Y\rangle$ is $C$.
  Then,
 conditions 1. and  2. in Theorem ~\ref{cpve} 
 are easily verified with $V$ the $C$-span of the $\X$ and $G=\GL$. Therefore,
 our main goal henceforth will be to show that 
 $F\langle \Y\rangle (\X)\supset F\langle \Y\rangle$ is a  no-new-constant extesion.
  
\subsection{The differential fields
 $\mathbf {F\langle \Y\rangle (\X)}$ and
 $\mathbf {F(\X)}$}\label{diffields}  

Let $\Y$, $1\leq i,j\leq n$,
 be differential indeterminates over $F$. For convenience, denote
 the $k$-th derivative 
$D^{(k)}(\Y)$ by $Y_{ij,k}$, for $k\geq 0$, so that 
 $D(Y_{ij,k})=Y_{ij,k+1}$,
  $k\geq 0$. As usual, $D^{(0)}(\Y)=Y_{ij,0}$ represents the
  original element $\Y$. In this situation  
we will omit the $k$-subindex  
and   write $\Y$ instead of $Y_{ij,0}$.

 Let $F\{\Y\}$ be the ring of differential polynomials in the $\Y$ and
 $F\langle\Y\rangle$ its differential quotient field. By that we mean
 the usual quotient field endowed with the natural derivation:
 $$\mathcal{D}\Big(  \frac pq\Big)= \frac{D(p)q-pD(q)}{q^2}.$$
 for $p, q\in F\{\Y\}$, where $D$ is the derivation on $F\{\Y\}$.
 
  Next let $\X$, $1\leq i,j\leq n$, be 
  algebraically independent over $F\langle \Y\rangle$.
  We consider the differential ring $R=F\{\Y\}[\X]$  with  
derivation extending the derivation on $F\{\Y\}$ by a formula
 $$D(\X)=\sum_{\ell=1}^nY_{i\ell}X_{{\ell}j}.$$
 
  As above, this
derivation extends to the quotient field  
$$Q=F\langle \Y\rangle(\X)$$  in  a natural way.

 Henceforth we  fix the differential field $Q$ as defined above. 

 Likewise we will regard the polynomial ring $F[\X]$ as a differential
 ring  with derivation extending the derivation on $F$ by a
 formula 
 $$D(\X)=\sum_{\ell=1}^nf_{i\ell}X_{{\ell}j}$$
with $\f\in F.$ Then, as before, we extend this
 derivation to the quotient field $F(\X)$ in a natural way.

The
multinomial notation 
$a_{\ba}\bx^{\ba}$ will be used
to denote a term of the form\linebreak
$a_{\alpha_{11}\cdots \alpha_{k\ell}}X^{\alpha_{11}}_{11}\cdots
X_{k\ell}^{\alpha_{k\ell}}.$

 The  ring $F[\X]$ is assumed to be ordered with  the degree 
reverse lexicographical  order (\emph{degrevlex}). That is,  the set 
$$\mathbb T^{n^2}=\lbrace \bx^{\bb}\,| \,
 \bx=(\X),\,\bb=(\beta_{ij})\in \mathbb N^{n^2}\rbrace$$ 
of the power products in the $\X$ is ordered by 
 $X_{11}>\dots>X_{1n}>\dots>
X_{n1}>\dots >X_{nn},$ and 
$$\bx^{\ba}< \bx^{\bb}\Longleftrightarrow\cases 
\sum_{j=1}^n\sum_{i=1}^n
\alpha_{ij}<\sum_{j=1}^n\sum_{i=1}^n \beta_{ij}\\
\text{\it or} \\
 \sum_{j=1}^n\sum_{i=1}^n \alpha_{ij}=\sum_{j=1}^n\sum_{i=1}^n 
\beta_{ij}, \text{\it and  the first coordinates}\\
\text{\it  $\alpha_{ij}$, $\beta_{ij}$ from the 
right which are different satisfy $\alpha_{ij}>\beta_{ij}$.}\endcases$$
Henceforth the leading power product of a polynomial in $F[\X]$  
is assumed to be with respect to this order.

The above definitions as well as the \emph{Multivariable Division Algorithm} that will
be used in Remark ~\ref{rmk4}, can be found in
\cite{A-L}.

\subsection{ Darboux polynomials and the 
 constants of 
 $\mathbf {F\langle \Y\rangle (\X)}$}  

We need to 
show that the field of constants $\C$ of 
 $Q=F\langle \Y\rangle (\X)$ coincides with
the field of constants $C$ of $F$. In this section we will show (Corollary ~\ref{maincor}) that
this can be reduced to proving that the only \begin{emph}{ Darboux polynomials}\end{emph}
in $R$ are, up to a scalar multiple in $F$,
powers of $\det[\X]$.

\begin{defi}\label{dp} Let $D$ be a
derivation on the polynomial ring $A=k[Y_1,$ $\dots,Y_s]$.  A polyomial 
$p\in A$  is called a Darboux polynomial if
  there is a  polynomial $q\in A$ such that
$D(p)=qp$. That is, $p$ divides $D(p)$.\end{defi}

Darboux polynomials correspond to generators of principal differential
ideals in $A$.  Chapter I of J.A. Weil's Thesis \cite{W} is devoted to
constants and Darboux polynomials in Differential Algebra. 

The following basic proposition (proven in \cite{W} for $A$ 
as in Definition ~\ref{dp}) characterizes    new constants
 for the extension $Q\supset F$ in terms of Darboux
polynomials:

\begin{proposition}\label{weil2}
 Let $p_1,p_2\in R=F\{\Y\}[\X]$, $p_1,p_2\ne 0$,
be relatively prime. Then $D(\frac {p_1}{p_2})=0$,
if and only if $p_1$ and $p_2$ are Darboux polynomials. Moreover, if $q_1,q_2\in R$
 are such that $D(p_1)=q_1p_1$ and $D(p_2)=q_2p_2$, then $q_1=q_2$. 
\end{proposition}

\begin{proof}
For the necessity of the condition we have 
$$D\left(  \frac {p_1}{p_2}\right)= \frac{D(p_1)p_2-p_1D(p_2)}{p_2^2} =0,$$
thus $D(p_1)p_2-p_1D(p_2)=0$, that is 
\begin{equation}D(p_1)p_2=p_1D(p_2). \tag{1}\end{equation}
Since $p_1$ and $p_2$ are
relatively prime, the last equation implies that $p_1$ divides $D(p_1)$ and
$p_2$ divides $D(p_2)$. 

Now, let $q_1,q_2\in R$ be such that $D(p_1)=q_1p_1$ and $D(p_2)=q_2p_2$, 
respectively. Then it follows from (1) that  
$$  q_1p_1p_2=q_2p_1p_2.$$
Hence, $q_1=q_2$.

The proof of the converse is obvious.
\end{proof}

\begin{proposition}\label{detdp}
Let $p=\ell\det[\X]^a$, with $\ell\in F$.
Then $p$ is a Darboux polynomial in $R$ with $q=\frac {\ell'}{\ell}+a\,\sum_{i=1}^n Y_{ii}.$
\end{proposition}
\begin{proof}
 We have
$$D(p)=\ell'\det[\X]^a+\ell( a\det[\X]^{a-1}D(\det[\X])).$$
And, $$D(\det[\X])=\Big(\sum_{i=1}^n Y_{ii}\Big) \det[\X].$$ 
Thus, 
\begin{eqnarray*} D(p)&=&\ell'\det[\X]^a+
\ell( a\det[\X]^{a-1}\Big(\sum_{i=1}^n Y_{ii}\Big)\det[\X]\\
&=&\Big(\frac {\ell'}{\ell}+a\,\sum_{i=1}^n Y_{ii}\Big)\,\ell\,\det[\X]^a.
\end{eqnarray*}
That is, $$p=\ell\det[\X]^a$$ 
 is a Darboux polynomial in $R$ 
 with $$q=\frac {\ell'}{\ell}+a\,\sum_{i=1}^n Y_{ii}.$$\end{proof}

\begin{corollary}\label{maincor}Suppose that if $p\in R$ is a Darboux polynomial then
$p=\ell\det[\X]^a$, with $\ell\in F$, $a\in\mathbb N$ and $q=\frac {\ell'}{\ell}+a\,\sum_{i=1}^n Y_{ii}.$
Then $F\langle\Y\rangle(\X)\supset F$ is a no-new-constant extension.
\end{corollary}

\begin{proof} $F\langle\Y\rangle(\X)$ is the fraction field of $F\{\Y\}[\X]$. Thus, if an 
element $f\in F\langle\Y\rangle(\X)$ satisfies $D(f)=0$ and $f\notin F$ then, by proposition ~\ref{weil2},
 $f= {p_1}/{p_2}$ 
with $p_1,p_2\in R$ relatively prime  Darboux polynomials. If the hypothesis is true, this 
contradicts the fact that $p_1$
and $p_2$ are relatively prime.\end{proof} 

Next we   show that the  constants of $F\langle\Y\rangle$ 
coincide with the  constants of $F$.  For simplicity, if $h(Y)\in F\{\Y\}$, we
will use the notation $h'(Y)$ for  $D(h(Y))$. Notice
that this is not  the usual meaning $h'(Y)=\sum h'_\alpha\mathbf Y^\alpha$.

\begin{proposition}\label{nodpy} If $h(Y)\in F\{\Y\}$
satisfies $h'(Y)=g(Y)h(Y)$ for some $g(Y)\in
F\{\Y\}$ then $h(Y)\in F$. That is, there are no non-trivial
Darboux polynomials in $F\{\Y\}$.\end{proposition} 

\begin{proof} 
Suppose that the hypothesis of the proposition is true.
According to our  notation  $D^{(k)}(\Y)=Y_{ij,k}$. Consider the set of subindices $\{ij,k\}$,
$i,j,k\in\mathbb N$,  
 ordered with the lexicographical order.
That is,  $\{i_1 j_1,k_1\}>\{i_2j_2,k_2\}$ if and only if the first
coordinates
$s_1$ and $s_2$ from the left,  for $s=i,j,k$ above, which are different satisfy
$s_1>s_2$. 
 
Let $\{mn,t\}$ be the largest subindex such that $Y_{mn,t}$ occurs in
$h(Y)$. 

Let $h(Y)= \sum_\alpha
a_{\alpha}Y^{\alpha_{11}}_{11}\cdots Y_{mn,t}^{\alpha_{mn,t}}.$
Then \begin{eqnarray*} h'(Y)&=& \sum_\alpha
a'_{\alpha}Y^{\alpha_{11}}_{11}\cdots Y_{mn,t}^{\alpha_{mn,t}}
+ \sum_\alpha
a_{\alpha}\alpha_{11}Y^{\alpha_{11}-
1}_{11}Y^{\alpha_{11,1}+1}_{11,1}\cdots Y_{mn,t}^{\alpha_{mn,t}}\\
& &\qquad\qquad\qquad\qquad+\dots + \sum_\alpha
a_{\alpha}\alpha_{mn,t}Y^{\alpha_{11}}_{11}\cdots 
Y_{mn,t}^{\alpha_{mn,t}-1}Y_{mn,t+1}\\
&=&h_1(Y_{11},\cdots ,Y_{mn,t})+\Big( \sum_\alpha
a_{\alpha}\alpha_{mn,t}Y^{\alpha_{11}}_{11}\cdots
 Y_{mn,t}^{\alpha_{mn,t}-
1}\Big)Y_{mn,t+1}\\
&=&g(Y)h(Y).\end{eqnarray*}

For $Y_{mn,t+1}=D(Y_{mn,t})$ we have $\{mn,t+1\}>\{mn,t\}$, thus it 
does not occur in $h(Y)$ by the choice of $\{mn,t\}$. Also, it does not
occur in $h_1(Y_{11},\cdots,$ $Y_{mn,t})$. Thus the above equation  
implies that $Y_{mn,t+1}$  must occur in $g(Y)$. Let $g_{t+1}(Y)$ be its
coefficient in $g(Y)$ and let 
$$h_2(Y)= \sum_\alpha
a_{\alpha}\alpha_{mn,t}Y^{\alpha_{11}}_{11}\cdots
 Y_{mn,t}^{\alpha_{mn,t}-1}.$$ Then we have
$$h(Y)g_{t+1}(Y)Y_{mn,t+1}=h_2(Y)Y_{mn,t+1}$$
or
$$h(Y)g_{t+1}(Y)=h_2(Y).$$

But the total degree of $h_2(Y)$ is strictly less than the total degree of $h(Y)$. 
This forces $h(Y)\in F.$\end{proof} 

\begin{corollary}\label{constfy}The field of constants of $F\langle\Y\rangle$ coincides
with $C$, the
field of constants of $F$.\end{corollary}
\begin{proof}
This is a consequence of Propositions ~\ref{weil2} and ~\ref{nodpy}.\end{proof}

\subsection{Darboux polynomials in $\mathbf{R=F\{\Y\}[\X]}$}
 Proposition ~\ref{detdp}  shows 
that  scalar multiples of $\det[\X]$ and its powers are Darboux polynomials in $R$.
Corollary ~\ref{maincor} implies that if
these are the only Darboux polynomials in $R$ then we are done since consequently there
will be no new constants in $Q$. In this section we will show that that is the case, namely,
the only Darboux polynomials in $R$ are those of the form 
 $p=\ell\det[\X]^a $, with $\ell\in F$, $a\in\mathbb N$ and
 $q=\frac {\ell'}{\ell}+a\,\sum_{i=1}^n Y_{ii}.$

\begin{remarks}
\begin{rmk}[Derivative 
of a power product in the $\X$]\label{rmk1} Let
 $$\bx^\alpha=
 X_{11}^{\alpha_{11}}\cdots X_{1n}^{\alpha_{1n}}\cdots X_{n1}^{\alpha_{n1}}\cdots X_{nn}^{\alpha_{nn}},$$
then 
\begin{eqnarray*}
\lefteqn{D(\bx^\alpha)=
\Big(\sum_{i=1}^n\sum_{j=1}^n \alpha_{ij}Y_{ii}\Big)\,
\bx^\alpha}\\
&& +  \sum_{i=1}^n\sum_{j=1}^n\Big(\sum_{\ell>i}\alpha_{ij}Y_{i\ell}X_{11}^{a_{11}}\cdots
\X^{\alpha_{ij}-1}\cdots X_{\ell j}^{\alpha_{\ell j}+1}\cdots X_{nn}^{\alpha_{nn}}\\
& &\qquad\qquad+\sum_{\ell<i}\alpha_{ij}Y_{i\ell}X_{11}^{a_{11}}\cdots
X_{\ell j}^{\alpha_{\ell j}+1}\cdots
\X^{\alpha_{ij}-1}\cdots X_{nn}^{\alpha_{nn}}\Big).
\end{eqnarray*}\end{rmk}

\begin{rmk}\label{rmk2} Given  $\bx^{\alpha}$ as in ~\ref{rmk1}, we need to know which 
power products in the $\X$ contain in their derivatives   a $Y$-multiple of $\bx^{\alpha}$.
That is, we need to find the power products $\bx^\beta$ such that $D(\bx^\beta)$ 
contains an expression of the form
$Y_{rt}\bx^{\alpha}.$
By  Remark ~\ref{rmk1} such power products are 
$$ \bx^{\alpha_{rs,t}}=\cases X_{11}^{\alpha_{11}}\cdots X_{rs}^{\alpha_{rs}+1}\cdots 
X_{ts}^{\alpha_{ts}-1}\cdots X_{nn}^{\alpha_{nn}}\,\,\,&\text{if $\,\,\,r<t$}\\
X_{11}^{\alpha_{11}}\cdots X_{ts}^{\alpha_{ts}-1}\cdots X_{rs}^{\alpha_{rs}+1}\cdots 
 X_{nn}^{\alpha_{nn}}\,\,\,&\text{if $\,\,\,r>t$}\endcases$$
for $1\leq r,s\leq n$, $t\ne r$, and $\bx^{\alpha}$ itself.\end{rmk}

\begin{rmk}\label{rmk3} Let $p\in R.$ 
Since  $D(\X)=\sum_{\ell=1}^nY_{i\ell}X_{\ell j}$, then the total
degree of $p$ with respect to the $\X$ does not change after
differentiation. Therefore, if 
$D(p)=qp$ then  
$q\in F\{\Y\}.$\end{rmk}\end{remarks}

\begin{proposition}\label{coeffmon} Let $p\in R$. Write it as $p=\sum_\alpha
p_\alpha(Y)\bx^\alpha$, with $p_\alpha(Y)\in F\{\Y\}$.
Then for any $\alpha$  with
$p_\alpha(Y)\ne 0$,  the coefficient of $\bx^\alpha$ in
$D(p)$ is $$p'_{\alpha}(Y)+
p_{\alpha}(Y)\sum_{i=1}^n\sum_{j=1}^n \alpha_{ij}Y_{ii} 
+\sum_{i=1}^n\sum_{j=1}^n(\alpha_{ij}+1)\sum _{\ell\ne
i}p_{\alpha_{ij,\ell}}(Y)Y_{i\ell},$$ 
where $\alpha_{ij,\ell}$ is the exponent vector of the power product
$$\bx^{\alpha_{ij,\ell}}=\cases X_{11}^{\alpha_{11}}\cdots \X^{\alpha_{ij}+1}\cdots 
X_{\ell j}^{\alpha_{\ell j}-1}\cdots X_{nn}^{\alpha_{nn}} \,\,\,&\text{if $\,\,\,i<\ell$}\\
X_{11}^{\alpha_{11}}\cdots X_{\ell j}^{\alpha_{\ell j}-1}\cdots \X^{\alpha_{ij}+1}\cdots 
 X_{nn}^{\alpha_{nn}}\,\,\,&\text{if $\,\,\,\ell>i$}\endcases$$
as in Remark ~\ref{rmk2}.\end{proposition}

\begin{proof}
This is a direct consequence of Remarks ~\ref{rmk1}  and ~\ref{rmk2}.\end{proof}

\begin{proposition}\label{pinx} Let $p\in
R$ and suppose that
$D(p)=qp$, for some $q\in
F\{\Y\}$. Then
$p\in F[\X].$\end{proposition}
\begin{proof} 
Let $p=\sum_\alpha
p_{\alpha}(Y)\bx^{\alpha}.$
Then
\begin{eqnarray*}
D(p)&=&\sum_\alpha p'_{\alpha}(Y)
\bx^{\alpha}
+p_{\alpha}(Y)
D(\bx^{\alpha})\\
&=&qp\\
&=&\sum_\alpha q(Y)p_{\alpha}(Y)
\bx^{\alpha}.\end{eqnarray*}

By Proposition ~\ref{coeffmon}, for each $\alpha$ with $p_\alpha(Y)\neq 0$
the corresponding coefficient of $\bx^\alpha$ in $D(p)$ is 
\begin{eqnarray*} D(p)_\alpha=p'_{\alpha}(Y)&+&
p_{\alpha}(Y) \sum_{i=1}^n\sum_{j=1}^n \alpha_{ij}Y_{ii}\\ 
&+& \sum_{i=1}^n\sum_{j=1}^n(\alpha_{ij}+1)
\sum _{\ell\ne i}p_{\alpha_{ij,\ell}}(Y)Y_{i\ell}.\end{eqnarray*}

 Since   $D(p)=qp$, it must be 
$D(p)_\alpha=q(Y)p_\alpha(Y)$ or, equivalently, 
\begin{eqnarray*} q(Y)p_{\alpha}(Y)=p'_{\alpha}(Y)&+&p_{\alpha}
(Y) \sum_{i=1}^n\sum_{j=1}^n \alpha_{ij}Y_{ii}\\ 
&+& \sum_{i=1}^n\sum_{j=1}^n(\alpha_{ij}+1)\sum _{\ell\ne
i}p_{\alpha_{ij,\ell}}(Y)Y_{i\ell}.\end{eqnarray*}
This means that for each $\alpha$, the coefficient $p_{\alpha}(Y)$ of 
$\bx^\alpha$ in $p$ divides the expression
$$  p'_{\alpha}(Y)+\sum_{i=1}^n\sum_{j=1}^n(\alpha_{ij}+1)
\sum _{\ell\ne i}p_{\alpha_{ij,\ell}}(Y)Y_{i\ell}.$$
Thus, for each $\alpha$,  there is $u_{\alpha}(Y)$ such that 
$$p_{\alpha}(Y)u_{\alpha}(Y)= p'_{\alpha}(Y)+ \sum_{i=1}^n\sum_{j=1}^n(\alpha_{ij}+1)
\sum _{\ell\ne i}p_{\alpha_{ij,\ell}}(Y)Y_{i\ell}.$$ 

As in the proof of Proposition ~\ref{nodpy}, order
the triples $\{ij,k\}$, $i,j,k\in\mathbb N$, with the lexicographical 
order. Let $\{mn,t\}$ be the largest subindex
such that $Y_{mn,t}$ occurs in $p$. We have $D(Y_{mn,t})=Y_{mn,t+1}$
and $\{mn,t+1\}>\{mn,t\}$.

Now, for each $\alpha$ 
such that 
$Y_{mn,t}$ occurs in $p_{\alpha}(Y)$ we have that
$Y_{mn,t+1}$ will occur in $p'_{\alpha}(Y)$ but not in $p_\alpha(Y)$
 or in 
$$ \sum_{i=1}^n\sum_{j=1}^n(\alpha_{ij}+1)\sum _{\ell\ne
i}p_{\alpha_{ij,\ell}}(Y)Y_{i\ell}$$ by the choice of $\{mn,t\}$.
Therefore, it must occur in 
$p_\alpha(Y)u_\alpha(Y)$. 
  Let  
$$p_{\alpha}(Y)=\sum
a_{\beta}Y^{\beta_{11}}_{11}
Y^{\beta_{12}}_{12}\cdots Y^{\beta_{mn,t}}_{mn,t}$$
then \begin{eqnarray*}\lefteqn{p'_{\alpha}(Y) =\sum\,\,  
a'_{\beta}Y^{\beta_{11}}_{11}\cdots Y^{\beta_{mn,t}}_{mn,t}}\\
&&\qquad\quad+\sum a_{\beta}\,\,{\beta_{11}}\,\,Y^{\beta_{11}-
1}_{11}Y^{\beta_{11,1}+1}_{11,1}\cdots Y^{\beta_{mn,t}}_{mn,t}+\dots\\
 & &\qquad\quad\quad+\sum a_{\beta}\,\,{\beta_{mn,t}}\,\,Y^{\beta_{11}}_{11}\cdots Y^{\beta_{mn,t}-
1}_{mn,t}Y_{mn,t+1}.\end{eqnarray*}

So $Y_{mn,t+1}$ occurs in $p'_{\alpha}(Y)$ only in 
\begin{eqnarray*}\lefteqn{ \sum \,\,a_{\beta}\,\,{\beta_{mn,t}}\,\,Y^{\beta_{11}}_{11}\cdots 
Y^{\beta_{mn,t}- 1}_{mn,t}Y_{mn,t+1}}\\
&&=\Big(\sum a_{\beta}\,\,{\beta_{mn,t}}\,\,Y^{\beta_{11}}_{11}\cdots Y^{\beta_{mn,t}-
1}_{mn,t}\Big)Y_{mn,t+1}\\
&&=v(Y)Y_{mn,t+1}.\end{eqnarray*}

Since $Y_{mn,t+1}$ occurs in $p_\alpha(Y)u_\alpha(Y)$ and not in
$p_\alpha(Y)$ it must occur in $u_\alpha(Y)$. Let $u_{\alpha,t+1}(Y)$
be the coefficient of $Y_{mn,t+1}$ in $u_\alpha(Y)$. Then it has to
be $$p_\alpha(Y)u_{\alpha,t+1}(Y)Y_{mn,t+1}=v(Y)Y_{mn,t+1}.$$

The above equation implies that
$p_\alpha(Y)$ divides  $v(Y)$. But this is impossible since
the total degree of $v(Y)$  is strictly less than the total
degree of $p_{\alpha}(Y)$. This contradiction yields the result.\end{proof}  

\begin{lemma}\label{coeffy}    Let $p\in F[\X]$ and suppose that
there is  $q \in F\{\Y\}$
such that $D(p)=qp.$ Then
 $q$ is a linear polynomial in the $\Y$. If $\beta=(\beta_{ij})$ is sucht that
 $\bx^{\beta}$ occurs in $p$, then  for $1\leq i\leq n$ the 
coefficient of $Y_{ii}$  in $q$
is $ \sum_{j=1}^n \beta_{ij}$. In particular,
the sums $\sum_{j=1}^n \beta_{ij}$, for $1\leq i\leq n$, 
are independent
of the choice of $\bx^\beta$.\end{lemma}

\begin{proof}
 We have $p=\sum a_{\beta}\bx^{\beta},$ with $a_{\beta}\in F.$

Thus, \begin{eqnarray*} D(p)&=&\sum a'_{\beta}
\bx^{\beta}
+a_{\beta}D(\bx^{\beta})\\
&=&qp\\
&=&\sum q(Y)a_{\beta}
\bx^{\beta}.
\end{eqnarray*}

By Proposition ~\ref{coeffmon}, the coefficient of $\bx^\beta$ in $D(p)$ is 
$$a'_{\beta}+
a_{\beta}  \sum_{i=1}^n\sum_{j=1}^n \beta_{ij}Y_{ii}+
 \sum_{i=1}^n\sum_{j=1}^n(\beta_{ij}+1)\sum _{\ell\ne
i}a_{\beta_{ij,\ell}}Y_{i\ell}.$$

Hence, it must be  
$$ q(Y)a_{\beta}=a'_{\beta}+
a_{\beta}\Big(  \sum_{i=1}^n\sum_{j=1}^n \beta_{ij}Y_{ii}
+ \sum_{i=1}^n\sum_{j=1}^n(\beta_{ij}+1)\sum _{\ell\ne
i}a_{\beta_{ij,\ell}}Y_{i\ell}\Big).
$$

From this, 
  $$ q(Y)={\frac{a'_{\beta}}{a_{\beta}}}+
  \sum_{i=1}^n\sum_{j=1}^n \beta_{ij}Y_{ii}+
 \sum_{i=1}^n\sum_{j=1}^n(\beta_{ij}+1)\sum _{\ell\ne
i}{\frac{a_{\beta_{ij,\ell}}}{a_{\beta}}}Y_{i\ell}.
 $$

The coefficient of $Y_{ii}$  in the above
expression is  $  \sum_{j=1}^n\beta_{ij}$,  
for $1\leq i\leq n$. Since this  expression for $q$
is valid for any index $\beta$, the ``in particular" part follows immediately.\end{proof}

\begin{corollary}\label{corcoeffy} Let $p$ be as in Lemma ~\ref{coeffy}. Let $\bx^\alpha$
be the leading power product of $p$. Let $\bx^{\beta}$ 
be any power product with non-zero coefficient in $p$. Then 
$  \sum_{j=1}^n \beta_{ij}=\sum_{j=1}^n \alpha_{ij}$, for
$1\leq i\leq n$. Thus $p$ is homogeneous of degree 
$  \sum_{j=1}^n\sum_{i=1}^n \alpha_{ij}.$
\end{corollary}
\begin{proof} This is an immediate consequence of the ``in particular"
part in Lemma ~\ref{coeffy}.\end{proof}
\begin{corollary}\label{corcoeffy2} Let $p\in
F[\X]$ and suppose that
$D(p)=qp$, for some $q\in
F\{\Y\}$. Let
$\bx^\alpha$ be the leading power product of 
$p$, and let $\ell\in F$ be its coefficient. 
Then
$$q={\frac{\ell'}{\ell}}+
 \sum_{i=1}^n\sum_{j=1}^n\alpha_{ij}Y_{ii}. $$
\end{corollary}

\begin{proof}
By Proposition ~\ref{coeffmon} and since $D(p)=qp$,
 the coefficient of $\bx^\alpha$ in $D(p)$ is
\begin{equation}
\ell q=\ell'+\ell\,\Big( 
\sum_{i=1}^n\sum_{j=1}^n\alpha_{ij}Y_{ii}+ 
\sum_{i=1}^n\sum_{j=1}^n(\alpha_{ij}+1)\sum _{\ell\ne
i}p_{\alpha_{ij,\ell}}Y_{i\ell}\Big).
                  \tag{1}\end{equation}

The $p_{\alpha_{ij,k}}$ are the coefficients of the power products $\bx^{\alpha_{ij,k}}$
in$p$, with ${\alpha_{ij,k}}\neq\alpha$, such that $D(\bx^{\alpha_{ij,k}})$ contains an expression
of the form $Y_{st}\bx^\alpha$. 
By Remark ~\ref{rmk2}, these power products are 
$$\bx^{\alpha_{rs,t}}=\cases 
X_{11}^{\alpha_{11}}\cdots X_{rs}^{\alpha_{rs}+1}\cdots 
X_{ts}^{\alpha_{ts}-1}\cdots X_{nn}^{\alpha_{nn}}\,\,\,&\text{if $\,\,\,r<t$}\\
X_{11}^{\alpha_{11}}\cdots X_{ts}^{\alpha_{ts}-1}\cdots 
X_{rs}^{\alpha_{rs}+1}\cdots 
 X_{nn}^{\alpha_{nn}}\,\,\,&\text{if $\,\,\,r>t$}\endcases,$$
all of which violate Corollary ~\ref{corcoeffy} for $i=r$ and $i=t$. Therefore it
must be $p_{\alpha_{ij,k}}=0$, for all $1\leq i,j\leq n$;  $k\ne i$. 
But now, substituting back in (1), we see that  
$$ \ell q= \ell'+\ell  \sum_{i=1}^n\sum_{j=1}^n\alpha_{ij}Y_{ii}.$$
 Hence,  $$q=\frac{\ell'}{\ell}+  \sum_{i=1}^n\sum_{j=1}^n\alpha_{ij}
Y_{ii}.$$\end{proof} 

Our next step in order to show that the Darboux polynomials $p\in R$ 
 have the desired form  will  be to show that such a $p$ \begin{emph}{is not reduced with
 respect to}\end{emph}
 $\det[\X]$. For that we will show that the leading
power product of $p$ is a power of the leading power product of
$\det[\X]$. First, we have

\begin{lemma}\label{alphalemma} Let $p\in F[\X]$ be such that $D(p)=qp, q\in
F\{\Y\}$. Let $\bx^\alpha$ be its leading power product. Then 
$\alpha_{ij}=0$ for $j\ne
n-i+1$ and $\alpha_{i,n-i+1}>0$, $1\leq i\leq n$. That is,
$\bx^\alpha=X_{1n}^{\alpha_{1n}}X_{2,n-1}^{\alpha_{2,n-1}}
\cdots X_{n1}^{\alpha_{n1}}$.\end{lemma}
 
\begin{proof}
To prove that $\alpha_{ij}=0$ for $j\ne n-i+1$ we first show that
 $\alpha_{ij}=0$ for 
$j>n-k+1,$ $ i\geq k,$ $2\leq k\leq n$. Indeed, for $k=2$ we
 have $j>n-1$, so $j=n$ and 
$$D(\bx^\alpha)=\alpha_{nn} \sum_{k=1}^{n-1}
Y_{nk}X_{11}^{\alpha_{11}}\cdots X_{kn}^{\alpha_{kn}+1}\cdots 
X_{nn}^{\alpha_{nn}-1}\,+\dots$$

Since $q$ has no $\Y$ with $i\ne j$, each term in $D(\bx^\alpha)$
 containing
such a $\Y$ must be cancelled. In particular 
we need to cancel the terms containing 
 $$Y_{nj}X_{11}^{\alpha_{11}}\cdots X_{jn}^{\alpha_{jn}+1}\cdots
 X_{nn}^{\alpha_{nn}-1}$$
for $1\leq j\leq n-1$ above. For that we can only use the 
derivatives of power products of
the form  
\begin{eqnarray*}\lefteqn{\bx^{\alpha_{nl,j}}=}\\
&&X_{11}^{\alpha_{11}}\cdots X_{j\ell}^{\alpha_{j\ell}-1}\cdots
 X_{jn}^{\alpha_{jn}+1}\cdots X_{n1}^{\alpha_{n1}}\cdots X_{n\ell}^{\alpha_{n\ell}+1}\cdots X_{nn}^{\alpha_{nn}-1},
\ \ \ \ \ \ \ \ell<n.\end{eqnarray*} But these are all  strictly greater than $\bx^\alpha$ (the leading
power product of $p$),
and they may not occur in $p$. 
As a consequence, it has to be $\alpha_{nn}=0$. Now let $k>2$ be such that 
$\alpha_{in}=0$ for $i\geq k$. Then 
\begin{eqnarray*}\lefteqn{\bx^\alpha=}\\ 
&&X_{11}^{\alpha_{11}}\cdots X_{k-1,n}^{\alpha_{k-1,n}}
\cdots X_{k,n-1}^{\alpha_{k,n-1}}X_{k+1,1}^{\alpha_{k+1,1}}
\cdots X_{k+1,n-1}^{\alpha_{k+1,n-1}}\cdots 
X_{n,n-1}^{\alpha_{n,n-1}}\end{eqnarray*}
and
\begin{eqnarray*}\lefteqn{D(\bx^\alpha)= 
\alpha_{k-1,n}\Big( \sum_{i<k-1}Y_{k-1,i}X_{11}^{\alpha_{11}}\cdots 
X_{in}^{\alpha_{in}+1}
\cdots X_{k-1,n}^{\alpha_{k-1,n}-1}\cdots X_{n,n-1}^{\alpha_{n,n-1}}}\\
&+& \sum_{i>k-1}Y_{k-1,i}X_{11}^{\alpha_{11}}\cdots 
 X_{k-1,n}^{\alpha_{k-1,n}-1}\cdots X_{in}^{\alpha_{in}+1}\cdots
 X_{n,n-1}^{\alpha_{n,n-1}}\Big)+\dots \end{eqnarray*}

Likewise, we need to cancel all 
the terms in $D(\bx^\alpha)$ that contain  $Y_{k-1,i}$,
for $i\ne k-1$.  
In particular, we need to cancel 
$$Y_{k-1,i}X_{11}^{\alpha_{11}}\cdots X_{in}^{\alpha_{in}+1}
\cdots X_{k-1,n}^{\alpha_{k-1,n}-1}\cdots X_{n,n-1}^{\alpha_{n,n-1}},$$
for $i< k-1$.  For that  we
can only use the power products of the form 
\begin{eqnarray*}\lefteqn{\bx^{\alpha_{k-1,\ell,i}}=}\\
  &&X_{11}^{\alpha_{11}}\cdots X_{i\ell}^{\alpha_{i\ell}-1}
\cdots X_{in}^{\alpha_{in}+1}\cdots X_{k-1,\ell}^{\alpha_{k-1,\ell}+1}
\cdots X_{k-1,n}^{\alpha_{k-1,n}-1}\cdots X_{n,n-1}^{\alpha_{n,n-1}},\end{eqnarray*}
for $i<k-1$.

 But all of them are 
strictly greater than $\bx^\alpha$ and cannot occur in $p$. Thus,
it has to be $\alpha_{k-1,n}=0$. Since this argument
 is valid for any $k>2$, it
follows that $\alpha_{kn}=0$, for $2\leq k\leq n$.
 This makes the statement that 
$\alpha_{ij}=0$ for $j>n-k+1$, $i\geq k$, true for $k=2$.

Now assume that $k$ is such that $\alpha_{ij}=0$ for $j>n-k+1, i\geq k$. So 
\begin{eqnarray*}\lefteqn{\bx^\alpha =}\\ 
&& X_{11}^{\alpha_{11}}\cdots X_{1n}^{\alpha_{1n}}\cdots  
X_{k,n-k+1}^{\alpha_{k,n-k+1}} 
 X_{k+1,1}^{\alpha_{k+1,1}}\cdots X_{k+1,n-k+1}^{\alpha_{k+1,n-k+1}}\cdots
 X_{n,n-k+1}^{\alpha_{n,n-k+1}}\end{eqnarray*}
and for  $i> k$ 
$$\alpha_{i,n-k+1}\Y X_{11}^{\alpha_{11}}\cdots 
X_{1n}^{\alpha_{1n}}\cdots 
X_{k,n-k+1}^{\alpha_{k,n-k+1}+1}\cdots 
 X_{i,n-k+1}^{\alpha_{i,n-k+1}-1}\cdots X_{n,n-k+1}^{\alpha_{n,n-k+1}}$$ 
occurs in $D(\bx^\alpha)$. Thus we need to cancel it. For that 
we can only use the derivatives of power products of the form

\begin{eqnarray*}\lefteqn{\bx^{\alpha_{ij,k}}=}\\
&&X_{11}^{\alpha_{11}}\cdots X_{kj}^{\alpha_{kj}-1}\cdots 
X_{k,n-k+1}^{\alpha_{k,n-k+1}+1}\cdots \X^{\alpha_{ij}+1}  
\cdots
  X_{i,n-k+1}^{\alpha_{i,n-k+1}-1}\cdots X_{n,n-k+1}^{\alpha_{n,n-k+1}}\end{eqnarray*}
with $j<n-k+1$ since  $\alpha_{kj}=0$ for all $j>n-k+1$ by hypothesis. But
all such power products are strictly greater than $\bx^\alpha$
 and therefore
they cannot occur in $p$. This forces $\alpha_{i,n-k+1}=0$ for $i>k$. We can
repeat this process until $k=n$ and get  $\alpha_{ij}=0$ for all 
$j>n-k+1$, $i\geq k$, $2\leq k\leq n$, that is, 
\begin{eqnarray*}\lefteqn{\bx^\alpha=}\\
&& X_{11}^{\alpha_{11}}\cdots X_{1n}^{\alpha_{1n}}
 X_{21}^{\alpha_{21}}\cdots 
X_{2,n-1}^{\alpha_{2,n-1}}X_{31}^{\alpha_{31}}\cdots X_{n-1,2}^{\alpha_{n-1,2}}
   X_{n1}^{\alpha_{n1}}.\end{eqnarray*}
 
Now we show that $\alpha_{ij}=0$ for $j<n-k+1$, $1\leq k\leq n-1$, $i\leq k$. 
The process is analogous
to what we just did. First we show that $\alpha_{i1}=0$ for $i<n$. Indeed,
for each $i$ we have for $\ell>i$ that 
 $$\alpha_{i1}Y_{i\ell} X_{11}^{\alpha_{11}}\cdots
  X_{i1}^{\alpha_{i1}-1}\cdots 
X_{\ell 1}^{\alpha_{\ell 1}+1}\cdots X_{n1}^{\alpha_{n1}} $$
occurs in $D(\bx^\alpha)$. So, in order to cancel it, 
we need to use the derivatives of power products of the form 
\begin{eqnarray*}\lefteqn{\bx^{\alpha_{ij,\ell}}=}\\
&&X_{11}^{\alpha_{11}}\cdots 
 X_{i1}^{\alpha_{i1}-1}\cdots
 \X^{\alpha_{ij}+1}\cdots
X_{\ell 1}^{\alpha_{\ell 1}+1}\cdots  X_{\ell j}^{\alpha_{\ell j}-1}\cdots 
X_{n1}^{\alpha_{n1}}\end{eqnarray*}
with $j>1$, all of which are strictly greater than
 $\bx^\alpha$ if $\ell<n$, and
for $\ell=n$ we cannot simply have one of those since $\alpha_{nj}=0$
 for $j\ne 1$. 
Thus such power products cannot occur in $p$ and it has to be $\alpha_{i1}=0$ for
$i<n$.

 Let $k\leq n-1$ be such that 
$\alpha_{ij}=0$ for $j<n-k+1$, $i\leq k$. 
We have 
\begin{eqnarray*}\lefteqn{\bx^\alpha=}\\
&&X_{1,n-k+1}^{\alpha_{1,n-k+1}}\cdots
 X_{1n}^{\alpha_{1n}}\cdots X_{k,n-
k+1}^{\alpha_{k,n-k+1}}\cdots X_{n1}^{\alpha_{n1}}\end{eqnarray*}
and for all $i<k$, $\ell>i$, we have that
$$\alpha_{i,n-k+1}Y_{i\ell}X_{1,n-k+1}^{\alpha_{1,n-k+1}}\cdots
 X_{i,n-k+1}^{\alpha_{i,n-k+1}-1}
\cdots X_{\ell, n-
k+1}^{\alpha_{\ell, n-k+1}+1}\cdots X_{n1}^{\alpha_{n1}}$$
occurs in $D(\bx^\alpha)$
and in order to cancel it  we only have the derivatives
of power products of the form
\begin{eqnarray*}\lefteqn{\bx^{\alpha_{ij,\ell}}=} \\
 & &X_{1,n-k+1}^{\alpha_{1,n-k+1}}\cdots X_{i,n-k+1}^{\alpha_{i,n-k+1}-1}
\cdots \X^{\alpha_{ij}+1}
\cdots X_{\ell,n-k+1}^{\alpha_{\ell, n-k+1}+1}\cdots 
X_{\ell j}^{\alpha_{\ell j}-1}
\cdots X_{n1}^{\alpha_{n1}}\end{eqnarray*}
with $j>n-k+1$ since $\alpha_{ij}=0$ for $i\leq k$, $j<n-k+1$. 
 
For $\ell<k$, all these power products are strictly greater than $\bx^\alpha$
 and therefore
they cannot occur in $p$. For $\ell\geq k$ we cannot simply have
 such power products since
for  $\ell\geq k$, $\alpha_{\ell j}=0$ if $j>n-k+1$. Thus it has to be
 $\alpha_{i,n-k+1}=0$
for $i\leq k-1$.

We can repeat this process until $k=n-1$ and get  $\alpha_{ij}=0$,
 $j<n-k+1$, $i\leq k$, $1\leq k\leq n-1$. This completes the proof
of the first
part of the lemma. 

To prove that $\alpha_{i,n-i+1}\ne 0$, for all $1\leq i\leq n$,
suppose that there is $i$ such that $\alpha_{i,n-i+1}=0$ and let $j\ne i$ be such
that $\alpha_{j,n-j+1}\ne 0$. Then
$ D(\bx^\alpha)$
will contain $$\alpha_{j,n-j+1}Y_{ji}X_{1n}^{\alpha_{1n}}\cdots 
X_{j,n-j+1}^{\alpha_{j,n-j+1}-1}\cdots 
 X_{i,n-j+1}\cdots X_{n1}^{\alpha_{n1}}
+\dots\quad\quad\quad\quad\text{if $i>j$}$$
or 
$$\alpha_{j,n-j+1}Y_{ji}X_{1n}^{\alpha_{1n}}\cdots X_{i,n-j+1}\cdots
 X_{j,n-j+1}^{\alpha_{j,n-j+1}-1}\cdots
 X_{n1}^{\alpha_{n1}}
+\dots \quad\quad\quad\quad\text{if $i<j$.}$$
As noted above, since $q$ does not contain any $\Y$ with $i\ne j$,
we need to cancel  the terms in $D(p)$ involving either of the above. But 
that is impossible since $\alpha_{ij}=0$ for all $j$ and by Corollary ~\ref{corcoeffy}
all the power products 
$$X_{11}^{\beta_{11}}\cdots \X^{\beta_{ij}}\cdots X_{nn}^{\beta_{nn}}$$ 
in $p$ must have
$\beta_{ij}=0$ for $j=1,\dots,n$. In particular,
 we cannot have in $p$ power
products of the form
$\bx^{\alpha_{j,n-j+1,i}}$ as in Remark ~\ref{rmk2}\end{proof}
 
Next we show that the  exponents $\alpha_{st}$  of the $X_{st}$ in  $\bx^\alpha$, the
leading power product of $p$,
are all equal:

\begin{lemma}\label{lp(p)} Let $p\in F[\X]$ be such that $D(p)=qp$, $q\in
F\{\Y\}$. Let $$\bx^\alpha=
X_{1n}^{\alpha_{1n}}X_{2,n-1}^{\alpha_{2,n-1}}\cdots
X_{n1}^{\alpha_{n1}}$$ be its leading power product. 
Then $\alpha_{i,n-i+1}=\alpha_{1n}$,
for $i>1$, that is, if $a=\alpha_{1n}$, then 
$$\bx^\alpha=(X_{1n}X_{2,n-1}\cdots X_{n1})^a.$$
\end{lemma}
 
\begin{proof}
Let $\ell$ be the coefficient of $\bx^\alpha$ in $p$. 
We have
\begin{equation*}\begin{split}&\quad D(\ell X_{1n}^{\alpha_{1n}}X_{2,n-1}^{\alpha_{2,n-1}}\cdots
 X_{n1}^{\alpha_{n1}})
=\\
&\qquad\Big( \sum_{i=1}^n
 \alpha_{i,n-i+1}\ell Y_{ii}\Big)X_{1n}^{\alpha_{1n}}
X_{2,n-1}^{\alpha_{2,n-1}}\cdots
X_{n1}^{\alpha_{n1}}\\
&\quad\qquad 
 +\alpha_{1n}\ell \sum_{k\ne 1}Y_{1k}
X_{1n}^{\alpha_{1n}-1}\cdots X_{k,n-
k+1}^{\alpha_{k,n-k+1}}\cdots X_{kn}\cdots X_{n1}^{\alpha_{n1}}\\
 &+\ell \sum_{1<i}\alpha_{i,n-i+1}\sum_{k>i}\Y
X_{1n}^{\alpha_{1n}}\cdots 
X_{i,n-i+1}^{\alpha_{i,n-i+1}-1}\cdots X_{k,n-k+1}^{\alpha_{k,n-k+1}}
\cdots X_{k,n-i+1}\cdots
X_{n1}^{\alpha_{n1}}\\
 &+\ell \sum_{1<i}\alpha_{i,n-i+1}\sum_{k>i}
\Y X_{1n}^{\alpha_{1n}}
\cdots X_{k,n-i+1}\cdots
X_{k,n-k+1}^{\alpha_{k,n-k+1}}\cdots  X_{i,n-i+1}^{\alpha_{i,n-i+1}-1}\cdots
X_{n1}^{\alpha_{n1}}\\ 
 &+\ell'X_{1n}^{\alpha_{1n}}X_{2,n-1}^{\alpha_{2,n-1}}\cdots 
X_{n1}^{\alpha_{n1}}.\end{split}\end{equation*}
 
In order to cancel $$\alpha_{1n}\ell Y_{1k}
X_{1n}^{\alpha_{1n}-1}\cdots X_{k,n-
k+1}^{\alpha_{k,n-k+1}}\cdots X_{kn}\cdots X_{n1}^{\alpha_{n1}},
\ \ \ \ \ k\ne 1,$$
above, we can only use the derivatives of
the power product 
\begin{eqnarray*}\lefteqn{\bx^{\alpha_{1,n-k+1,k}}=}\\
&&X_{1,n-k+1}\cdots X_{1n}^{\alpha_{1n}-1}\cdots
 X_{k,n-k+1}^{\alpha_{k,n-k+1}-1}
\cdots X_{kn}\cdots X_{n1}^{\alpha_{n1}},\end{eqnarray*}
since for $j\ne n-k+1$ we have $\alpha_{kj}=0$.

 Let $a_{\alpha_{1,n-k+1,k}}$ be the coefficient of $X^{\alpha_{1,n-k+1,k}}$
 in $p$. Then 
\begin{equation}a_{\alpha_{1,n-k+1,k}}=-\ell \alpha_{1n} \tag{2}\end{equation}

On the other hand, in order to cancel 
$$\alpha_{k,n-k+1}\ell Y_{k1}X_{1,n-k+1}\cdots X_{1n}^{\alpha_{1n}}\cdots 
X_{k,n-k+1}^{\alpha_{k,n-k+1}-
1}\cdots X_{n1}^{\alpha_{n1}},\ \ \ \ \ k\ne 1$$
above, the only power product that we can use is, again,
\begin{eqnarray*} \bx^{\alpha_{kn,1}} 
&=&X_{1,n-k+1}\cdots X_{1n}^{\alpha_{1n}-1}\cdots
 X_{k,n-k+1}^{\alpha_{k,n-k+1}-1}\cdots X_{kn}\cdots
X_{n1}^{\alpha_{n1}}\\
&=& X^{\alpha_{1,n-k+1,k}},\end{eqnarray*}
since $\alpha_{1j}=0$ for $j\ne n$. Thus it must be 
\begin{equation}a_{\alpha_{1,n-k+1,k}}=-\ell \alpha_{k,n-k+1} \tag{3}\end{equation}
as well.

From (2) and (3) it follows that, for $k\ne
1$, $\alpha_{1n}=\alpha_{k,n-k+1}.$ \end{proof}
 
As a consequence of the above results we obtain the following expression
for $q$:

\begin{corollary}\label{corq}  Let $p\in
F[\X]$ and suppose that
$D(p)=qp$, $q\in
F\{\Y\}$. Let $\bx^\alpha$ be the
leading power product of $p$. Let $a\in\mathbb N$ be such that $$\bx^\alpha=(X_{1n}X_{2,n-1}\cdots 
X_{n1})^a$$
 and let $\ell\in F$ be the coefficient of $\bx^\alpha$ in $p$. 
Then
$$q=\frac{\ell'}{ \ell}+ a\, \sum_{i=1}^nY_{ii}. $$
\end{corollary}

\begin{proof}
This is a consequence of Corollary ~\ref{corcoeffy2} and Lemma ~\ref{lp(p)}.\end{proof}

\begin{corollary}\label{phomogen} Let $p$ be as in Corollary ~\ref{corq}. 
Then $p$ is homogeneous
of degree $na.$\end{corollary}

\begin{proof} 
This is a consequence of Corollary ~\ref{coeffy} and Lemma ~\ref{lp(p)}\end{proof}  

Lemma ~\ref{lp(p)} implies that $p$ is not reduced with respect
to $\det[\X]$. Since this is a key point in the
proof of our main result we restate it 
as the following

\begin{theorem}\label{bigthm}  Let $p\in F[\X]$ be such that $D(p)=qp$,
$q\in F\{\Y\}$. Let $\bx^\alpha$ be its leading power product. Then
$$\bx^\alpha=(X_{1n}X_{2,n-1}\cdots X_{n1})^a=
\text{\rm lp}(\det[\X])^a.$$ 
Thus 
$p$ is not reduced with respect to $\det[\X]$.
\end{theorem}

{\bf Note.} If $f$ is a polynomial, $\text{\rm lp}(f)$ denotes its leading
power product with respect to a given order.

\begin{proof}
This is just a restatement of Lemma ~\ref{lp(p)}.\end{proof}

\begin{remark}\label{rmk4}  Let $p_1,p_2\in F[\X]$ be two polynomials such that
$\text{lp}(p_1)=X^\alpha=\text{lp}(p_2)$. Then we can write $p_1=f\,p_2+r$ 
where $f\in F$ and
 $r$ is reduced with respect to $p_2$. Indeed, since
$\text{lp}(p_1)=\text{lp}(p_2)$, we have that $\text{lp}(p_2)$ divides 
$\text{lp}(p_1)$.
So $p_1$ is not reduced with respect to $p_2$ . We may apply the
Multivariable Division Algorithm to $p_1$ and $p_2$, 
to get $f,r\in F[\X]$, such that $p_1=f\,p_2+r$, with $r$
 reduced with respect to $p_2$ and 
$\text{lp}(p_1)=\text{lp}(f)\text{lp}(p_2)$. The last equation implies that
$\text{lp}(f)=1$. Hence, $f\in F$.\end{remark}

We are now ready to prove our main result
on the form of the Darboux polynomials in $R$:

\begin{theorem}\label{bgstthm}   Let $p\in F[\X]$ and $q\in F\{\Y\}$
be polynomials in $R$ that satisfy the Darboux condition $D(p)=qp$.
 Then there is $a\in \mathbb N$ and
 $\ell\in F$ such that  $$p=\ell\,\det[\X]^a$$ 
and 
$$q=\frac{\ell'}{ \ell}+a\, \sum_{i=1}^nY_{ii}. $$
\end{theorem}
 
\begin{proof}
Let $q_1=\sum_{i=1}^n a\,Y_{ii},$ so that,
$$ D(\det[\X]^a)=q_1\det[\X]^a
=(q-\frac{\ell'}{ \ell})\det[\X]^a.$$

By Remark ~\ref{rmk4} we can write $p=\ell\det[\X]^a+r$, with $r$ reduced
with respect to $\det[\X]^a$. Now, 
\begin{eqnarray*} D(p)&=&D(\ell\det[\X]^a)
+D(r)\\
&=&\ell'\det[\X]^a+\ell(q-\frac{\ell'}{\ell})\det[\X]^a+D(r)\\
&=&\ell'\det[\X]^a+q\ell\det[\X]^a-\ell'\det[\X]^a+D(r)\\
&=&q\ell\det[\X]^a+D(r).\end{eqnarray*}

On the other hand, we have 
\begin{eqnarray*}  D(p)&=&qp\\ &=&q\ell\det[\X]^a+qr.\end{eqnarray*}
 
Therefore, it has to be $D(r)=qr$.
But $r$ is reduced with respect to $\det[\X]^a$. It follows, by
Theorem ~\ref{bigthm}, that $r=0$.
The statement about the form of $q$ is just the
content of  Corollary ~\ref{corq}.\end{proof}
 
Now, for $1\leq i,j\leq n$, 
let $D_{E(ij)}\in \Lie(\GL)$ be the derivation given by multiplication by the matrix $E(ij)$, with 
1 in position (i,j) and zero elsewhere. The set $\{D_{E(ij)}|\,1\leq i,j\leq n\}$ is a basis for
$\Lie(\GL)$ with respect to which the derivation in Theorem ~\ref{bgstthm}  is expressed.
We show that the result  does not depend on the basis chosen on $\Lie(\GL)$:  

\begin{theorem}\label{thmingen} Let
${\mathcal D}_{st}$, $1\leq s,t\leq n$, be any  basis of $\Lie(\GL)$.
Define a derivation in the ring 
$R=F\{\Y\}[\X]$ 
by ${\mathcal D}= \sum Y_{st}{\mathcal D}_{st}$. Let $p$ and $q$ be 
 polynomials in $R$ that satisfy the Darboux condition
  ${\mathcal D}(p)=qp$. Then there is
 $a\in\mathbb N$ and  $\ell\in F$ such that
$p=\ell\det[\X]^a$ and 
$q=\frac{\ell'}{\ell}+a\, \sum_{i=1}^nY_{ii}$.\end{theorem}

\begin{proof}
Since  $\{D_{E(ij)}|\,1\leq i,j\leq n\}$  is a basis  of $\Lie(\GL)$ we have
$${\mathcal D}_{st}= \sum c_{st,ij}D_{E(ij)}, $$
with $c_{st,ij}\in C$. Thus, \begin{eqnarray*}  {\mathcal D}&=& \sum_{s,t}
Y_{st}{\mathcal D}_{st}\\
&=& \sum_{s,t} Y_{st} \sum_{i,j} c_{st,ij}D_{E(ij)}\\
&=& \sum_{i,j}\sum_{s,t}c_{st,ij}Y_{st}D_{E(ij)}\\
&=& \sum_{i,j}Z_{ij}D_{E(ij)},\end{eqnarray*}
where $Z_{ij}= \sum_{s,t}c_{st,ij}Y_{st}$. Now, 
$[c_{st,ij}]$ is a matrix of change of basis so it is invertible. Also
the 
$c_{st,ij}$
are contants for $D$, thus
the map $Z_{ij,k}\rightarrow Y_{ij,k}$ is a  differential bijection. 
In other words, the differential rings 
$$R=F\{\Y\}[\X], D$$ and
$$R'=F\{Z_{ij}\}[\X], {\mathcal D}$$ are 
isomorphic and therefore we can apply Theorem ~\ref{bgstthm} to $R'.$\end{proof}


\section{A Generic  Picard-Vessiot  Extension for $\GL$}

In this  section we prove the statement about the generic Picard-Vessiot  extension for $F$
with differential Galois group
 $\GL$.  We  provide some specialization properties 
  in the
 case when $F$ has finite transcendence degree over  $C$.

\subsection{The generic extension}  

\begin{theorem}\label{thmgpve} $F\langle \Y\rangle(\X)\supset F\langle \Y\rangle$
is a generic Picard-Vessiot  extension  with differential
Galois group $\GL$.
 \end{theorem}

\begin{proof}
First we
need to show that $F\langle \Y\rangle(\X)\supset F\langle \Y\rangle$ is a Picard-Vessiot
extension with differential Galois group $\GL$.  We will
use the characterization of Theorem ~\ref{cpve}. We have
\begin{enumerate}
\item $F\langle \Y\rangle(\X)=F\langle \Y\rangle\langle V\rangle$,
 where $V\subset F\langle \Y\rangle(\X)$
 is the finite dimensional
vector space over $C$ spanned by the $\X$.
\item The group $G=\GL$ acts as a group of differential
automorphisms of\linebreak $F\langle \Y\rangle(\X)$ with $G(V)\subseteq V$ and $F\langle \Y\rangle(\X)^G=F\langle
\Y\rangle$. This follows from the fact that $F\langle \Y\rangle(\X)$ is
the function field of $\GL_{F\langle \Y\rangle}$.

\item  $F\langle \Y\rangle(\X)\supseteq F\langle \Y\rangle$ has 
no new constants. This is a consequence of Proposition ~\ref{weil2}, Corollary ~\ref{maincor} 
and Theorem ~\ref{bgstthm}.\end{enumerate}

Now, suppose that  $E\supseteq F$ is 
a Picard-Vessiot  extension of $F$ with differential Galois
group $\GL$.
By Theorems ~\ref{sthm} and ~\ref{stgln},
 we have that in this situation $E$ is isomorphic to $F(\X)$ 
(the function field of $\GL_F$) as a $\GL$-module
 and as an $F$-module. Any $\GL$
equivariant 
  derivation $D_E$ on $F(\X)$ extends the derivation on $F$ in such a way that
 $$D_E(\X)=\sum_{\ell=1}^n f_{i\ell}X_{\ell j}$$ with
$\f\in F$. Since $E\supset F$ is a Picard-Vessiot  extension for
$\GL$, then
so is $C\langle \f\rangle(\X)\supset C\langle \f\rangle$,
 the derivation on  $C\langle \f\rangle(\X)$   being 
the corresponding restriction of $D_E$. From this Picard-Vessiot 
extension one can retrieve $F(\X)\supset F$ by extension of scalars from
$C$ to $F$.  In this way, any Picard-Vessiot extension
$E\supset F$ with 
differential Galois group $\GL$ can be obtained from 
$F\langle \Y\rangle(\X)\supset F\langle \Y\rangle$
via the specialization $\Y\mapsto \f$. This means that $F\langle \Y\rangle(\X)\supset F\langle \Y\rangle$
is a generic Picard-Vessiot extension of $F$ for $\GL$.\end{proof}
  
\subsection{Specializing to a Picard-Vessiot  extension of $\mathbf F$}\label{sectsp}
   
 In this section we give necessary and sufficient conditions
for a specialization $\Y\mapsto \f$, $\f\in F$, with
$C\langle \f\rangle(\X)\supset C\langle \f\rangle$
a Picard-Vessiot  extension, to exist. We restrict ourselves
to the case when $F$ has finite transcendence degree over
$C$. 

  Our goal is to find $f_{ij}\in F$ such that  
the specialization (homomorphism) from $C\{\Y\}$ to $F$ given 
by $\Y\mapsto \f$
is such that  $C\langle \f\rangle(\X)\supset C\langle \f\rangle$,
 with 
derivation given by $D(\X)= \sum_{\ell=1}^nf_{i\ell}X_{\ell j}$, has no new constants. We have:

\begin{theorem}\label{specializ2}   
Let $F=C(t_1,\dots,t_m)[z_1,\dots,z_k]$
where the $t_i$ are algebraically independent over
$C$ and the
$z_i$ are algebraic over $C(t_1,\dots,t_m)$.
Assume that the derivation on $F$ has field of constants $C$ and that
it extends to $F(\X)$ so that 
$D(f\otimes \X)= D(f)\otimes \X+f\otimes \sum_{\ell=1}^n
f_{i\ell}X_{\ell j}$ on $F\otimes C[\X]$.
 Let $\C$
be the field of constants of $F(\X)$. Then $\C=C$ if and only if the set of all the
$t_i$ and all the $\X$
are algebraically independent over $\C$.\end{theorem}

\begin{proof} (Sufficiency) Suppose that $\C$ properly contains $C$. Let $r$ be
the transcendence degree of $\C$ over $C$. Since
$C$ is algebraically closed, $r$  has to be at least one.
 
We have the tower of fields
$$C\subset\C\subset\C(\X)\subset F(\X)$$
where the transcendence degree of $\C\subset\C(\X)$ is $n^2$ and 
the transcendence degree of $C\subset F(\X)$ is $n^2+m$. Since $r\geq 1$
the transcendence degree $\ell$ of $\C(\X)\subset F(\X)$
 has to be $\ell<m$ and therefore  there 
is an algebraic relation among the $t_i$
over $\C(\X)$. 
Let $g(\X), f_i(\X)\in \C[\X]$,
 $g(\X)\not\equiv 0$,
be such that 
$$ t^{\delta_s}+
 \frac{f_{s-1}(\X)}{g(\X)}t^{\delta_{s-1}}+\dots +\frac
{f_0(\X)}{g(\X)}=0.$$

Then 
$$ g(\X)t^{\delta_s}+
 f_{s-1}(\X)t^{\delta_{s-1}}+\dots +
f_0(\X)=0.$$

Since the $f_i(\X)$ and $g(\X)$ are polynomials in the $\X$ with coefficients in
$\C$, the last equation gives an algebraic relation among the $t_i$ and the $\X$ over
$\C$.

For the necessity we only need to point out that by construction the set of all the $t_i$
and all the $\X$
are algebraically independent over $C$.\end{proof}

Now to check whether the set of all the $t_i$ and all the $\X$ are algebraically independent over 
$\C$, we let $\mathbb T_k$, $k\geq 1$, denote the set of monomials in both the $t_i$ and the $\X$
 of total degree less than or equal to $k$.
 Then the set of all the $t_i$ and all the $\X$ are algebraically independent over 
$\C$ if and only if, for each $k$, the set $\mathbb T_k$ 
is linearly independent over $\C$.
 
Fix a term order on the set $\mathbb T$ of all monomials in both the $t_i$ and the $\X$ and let $W_k$
 denote the wronskian 
of the set $\mathbb T_k$ relative to that order.
Then the above 
condition is equivalent to the fact that  $W_k\neq 0$ for $k\geq 1$.
Now go back to $C\{\Y\}[\X]$ and extend scalars from $C$ to $F$.
Let $W_k(\Y)$ be the Wronskian of 
$\mathbb T_k$ in $F\otimes C\{\Y\}[\X]$. 

 Then, the condition of Theorem ~\ref{specializ2} for
finding a specialization $\Y\mapsto \f$
so that $C\langle \f\rangle(\X)\supset C\langle \f\rangle$ has no new constants can be expressed
as follows:

\begin{theorem}\label{wronskiancond}  There is a specialization of the $\Y$ with no
new constants if and only if  there are $\f\in F$ such that all the wronskians $W_k(\Y)$, $k\geq 1$, map
to non-zero elements under $\Y\mapsto\f$.\end{theorem}

\section{Results for  
 connected linear algebraic groups}
  We do not know at present whether generic Picard-Vessiot 
  extensions exist for arbitrary connected linear algebraic groups.  However, the proofs of the 
  specialization 
  theorems in ~\ref{sectsp} can be easily generalized for such groups.
\subsection{Specialization results}
We point out that the proofs of Theorems  ~\ref{specializ2} and\linebreak 
 ~\ref{wronskiancond}  do not make any special
use of the fact that $G=\GL$. 

Let $F=C(t_1,\dots,t_m)[z_1,\dots,z_k]$
where the $t_i$ are algebraically independent over $C$ and the
$z_i$ are algebraic over $C(t_1,\dots,t_m)$.

Let $Y_1,\dots,Y_n$ be differential indeterminates over $F$ and
$X_1,\dots,X_n$ algebraically independent over $F\langle Y_i\rangle$.

We consider the group $G$ to be a connected linear algebraic group over $C$
with function field  $C(G)=C(X_i)$.

If $\{D_1,\dots,D_n\}$ is a basis for $\Lie(G)$, let 
$D_Y= \sum_{i=1}^n Y_i\,D_i $ be a $G$-equivariant derivation 
on $F\langle Y_i\rangle(X_i)$ and 
$D= \sum_{i=1}^n f_i\,D_i $, $f_i\in F$, a specialization of $D_Y$ to 
a $G$-equivariant derivation on $F(X_i)$. Let $\C$
be the field of constants of $F(X_i)$ for this derivation. We have,

\begin{theorem}\label{specializ3} Let $F$, $C$ and $\C$ be as above.  
Then $\C=C$ if and only if the set of all the $t_i$ and the $X_i$
are algebraically independent over $\C$.\end{theorem}

Now, fix an order in the
 set $\mathbb T$ of monomials  in both the $t_i$ and the $X_i$ and let 
$W_k(Y_i)$ be the wronskian (with respect to this order)
 of the monomials in both the $t_i$ and the $X_i$ of degree less than or equal to $k$ 
 computed in $F\otimes C\{Y_i\}[X_i]$.
Then,  

\begin{theorem}\label{wronskiancond1}  There is a specialization of the $Y_i$ with no
new constants if and only if  there are $f_i\in F$ such that all the wronskians
$W_k(Y_i)$, $k\geq 1$, map
to non-zero elements under $Y_i\mapsto f_i$.\end{theorem}

For the proofs of Theorems ~\ref{specializ3}  and  ~\ref{wronskiancond1}    we only need to replace
the $\X$ with $X_i$ , the $\Y$ with $Y_i$  and $n^2$ with
$n$ in the proofs of Theorems ~\ref{specializ2} and ~\ref{wronskiancond}.$\hfill\qed$
 
 Observe that the proofs of Theorems ~\ref{specializ3} and  ~\ref{wronskiancond1}    do
not  use  the fact that  $C(X_i)$ is the
function field of $G$. However, this hypothesis is used in the 
following theorem to show that  $F(X_i)\supset F$ is
a Picard-Vessiot  extension with group $G$.

Under the hypothesis and notation of Theorems ~\ref{specializ3} and  ~\ref{wronskiancond1}   
we have:

\begin{theorem}\label{specializ4}  $F(X_i)\supset F$ is a Picard-Vessiot 
extension with Galois group $G$ if and only if the set of all the $t_i$ and all the
 $X_{i}$ are algebraically independent over
the field of constants $\C$ of $F(X_i)$. \end{theorem}

\begin{proof}
First assume that  $F(X_i)\supset F$ 
is a Picard-Vessiot
extension. Then the field of constants $\C$ of $F(X_i)$ coincides with $C$.
 So we can apply Theorem ~\ref{specializ3} and get the result.

Conversely, 
if the set of all the $t_i$ and all the $X_{i}$ are algebraically independent over
 $\C$,
by Theorem ~\ref{specializ3}, $F(X_i)\supset F$ is a no new constant extension. 
On the other hand, $F(X_i)$ is obtained from $C(X_i)$ by the extension of scalars: 
\begin{eqnarray*}  F(X_i)&=&q.f.(F\otimes_CC(X_i))\\
&=&q.f.(F\otimes_CC[G])\end{eqnarray*}
where $C[G]$ is the coordinate ring of $G$ and
$G$ acts on $F\otimes_CC[G]$ fixing $F$. So, $ G\subseteq G(F(X_i)/F)$.
Counting dimensions we get that $ G=G(F(X_i)/F) $ since
$C(X_i)=C(G)$, the function field fo $G$.
Finally, $F(X_i)=F\langle V \rangle$,
 where $V$ is the finite-dimensional vector space
over $C$ spanned by the $X_i$. By Theorem ~\ref{cpve},  $F(X_i)\supset F$
is a Picard-Vessiot extension.
\end{proof}

Applying Theorems   ~\ref{wronskiancond1}    and ~\ref{specializ4} we also obtain: 

\begin{theorem}\label{specializ5}  There is a specialization of the $Y_i$ such
that $F(X_i)\supset F$ is a Picard-Vessiot  extension if and only if  
there are $f_i\in F$ such that all the $W_k(Y_i)$, $k\geq 1$,
 map to non-zero elements via $Y_i\mapsto f_i$. 
\end{theorem}

\section{ Computing new constants}
   
 Let $C$ be an algebraically closed field with trivial derivation.
 Let\linebreak $F=C(t_1,\dots,t_m)[z_1,\dots,z_k]$
where the $t_i$ are algebraically independent over $C$ and the
$z_i$ are algebraic over $C(t_1,\dots,t_m)$.
Assume that the derivation on $F$ has field of constants $C$ and that
it extends to $F(\X)$ so that 
$$D(f\otimes \X)= D(f)\otimes \X+f\otimes \sum_{\ell=1}^n
f_{i\ell}X_{\ell j}$$ on $F\otimes C[\X]$, for certain $\f\in F$. 
By Theorem ~\ref{specializ2},  if there is an algebraic 
relation among the set of all the $t_i$ and all the $\X$ over the field of constants $\C$ 
of $F(\X)$ then $\C$
 properly contains $C$.

In this section we will produce a new constant from such an algebraic
relation. In order to simplify the computations we  will assume $F=C$.
So, in particular, the coefficients $f_{ij}$ in the derivation above are constant.
In this situation, since the transcendence degree of $F$ over $C$ is zero, if
$\C\varsupsetneqq C$, the 
condition of Theorem ~\ref{specializ2} means that the $\X$ are algebraically
dependent over $\C$.

We will restrict ourselves to the case $n=2$ and  use a particular linear
dependence relation. 
 
 Extend
the derivation on $F$ to $F(X_{11},X_{12},X_{21},X_{22})$
by letting 
 $$D(\X)=\sum_{\ell=1}^2f_{i\ell}X_{\ell j},$$ 
where the $\f$ are such that
 the wronskian
$W_1=w(X_{11},X_{12},$ $X_{21},X_{22})= 0$. That is,
the $\X$ are linearly dependent over $\C$.
Furthermore, assume that the linear 
relation among the $\X$ is such that  
there are $\beta_{12},\beta_{21},\beta_{22}\in \C$ with
\begin{equation}X_{11}=\beta_{12}X_{12}+\beta_{21}X_{21}+\beta_{22}X_{22} \tag{1}\end{equation}
and that $X_{12},X_{21}$ and $X_{22}$ are linearly independent.
In order to simplify the 
computations we will also assume that $\det[\f]=0$.

 We want to find $a,b,c\in F$  such that $p=aX_{12}+bX_{21}+cX_{22}$ is a Darboux
polynomial in $F[\X]$, that is 
$D(aX_{12}+bX_{21}+cX_{22})=q(aX_{12}+bX_{21}+cX_{22})$
for certain $q\in F$. 

 We have,
\begin{eqnarray*}\lefteqn{D(a\,X_{12}+b\,X_{21}+c\,X_{22})}\\&\quad=&a(f_{11}X_{12}+f_{12}X_{22})
+b(f_{21}X_{11}+f_{22}X_{21})
+c(f_{21}X_{12}+f_{22}X_{22})\\
&\quad=&bf_{21}X_{11}+(af_{11}+cf_{21})X_{12}+bf_{22}X_{21}+(af_{12}+cf_{22})
X_{22}\\
&\quad=&bf_{21}(\beta_{12}X_{12}+\beta_{21}X_{21}+\beta_{22}X_{22})+
(af_{11}+cf_{21})X_{12} +bf_{22}X_{21}\\
& &\qquad+(af_{12}+cf_{22})X_{22}\\
&\quad=&(af_{11}+bf_{21}\beta_{12}+cf_{21})X_{12}+b(f_{22}+f_{21}\beta_{12})X_{21}\\
& &\qquad+(af_{12}+bf_{21}\beta_{22}+cf_{22})X_{22}\\
&\quad=&qaX_{12}+qbX_{21}+qcX_{22}.\end{eqnarray*}
Therefore,
\begin{equation*}\begin{split} 
 [a(f_{11}-q)+bf_{21}\beta_{12}+cf_{21}&]X_{12}+b(f_{22}+
f_{21}\beta_{12}-q)X_{21}\\ 
&+(af_{12}+bf_{21}\beta_{22}+c(f_{22}-q)X_{22}
=0. \end{split}\tag{2}\end{equation*}  

Since we are assuming that $X_{12}, X_{21}$ and $X_{22}$ are linearly
independent their coefficients in (2) must be equal to zero. 
So we have the following homogeneous linear system in $a,b,c$:
$$\begin{array}{r@{\:\:}c@{\:\:}r@{\:\:}c@{\:\:}r@{\:\;=\:\;}c}
(f_{11}-q)\,a&+& f_{21}\beta_{12}\,b&+& f_{21}\,c & 0\\
& &(f_{22}+f_{21}\beta_{12}-q)\,b& & & 0\\
f_{12}\,a&+& f_{21}\beta_{22}\,b&+& (f_{22}-q)\,c & 0
\end{array}$$ 

In order for the above system to have non-trivial solutions we need that 
$$\det\left[\begin{array}{*{3}{c@{\quad}}}
f_{11}-q&f_{21}\beta_{12}& f_{21} \\
0&f_{22}+f_{21}\beta_{12}-q&0\\
f_{12} &f_{21}\beta_{22} & f_{22}-q  
\end{array}\right] =0.$$ 

But,
\begin{eqnarray*}\lefteqn{\det\left[\begin{array}{*{3}{c@{\quad}}}
f_{11}-q&f_{21}\beta_{12}& f_{21} \\
0&f_{22}+f_{21}\beta_{12}-q&0 \\
f_{12} &f_{21}\beta_{22} & f_{22}-q  
\end{array}\right]}\\
&&=(f_{22}+f_{21}\beta_{12}-q)\det\left[\begin{array}{*{2}{c@{\quad}}}
f_{11}-q& f_{21} \\
f_{12} & f_{22}-q\end{array}\right] \end{eqnarray*}
\begin{eqnarray*} 
&=&(f_{22}+f_{21}\beta_{12}-q)(\det[\f]-\Big(\sum_{i=1}^2f_{ii}\Big)q+q^2)\\ &=&0.\end{eqnarray*} 
 
This gives either
\begin{equation}f_{22}+f_{21}\beta_{12}-q=0 \tag{3}\end{equation} or 
\begin{equation}\det[\f]-\Big(\sum_{i=1}^2f_{ii}\Big)q+q^2=0. \tag{4}\end{equation}
From (3)-(4) we get 
\begin{equation}q=f_{22}+f_{21}\beta_{12} \tag{5}\end{equation} or
\begin{equation}q=\frac{\sum_{i=1}^2f_{ii}\pm\sqrt{\big(\sum_{i=1}^2f_{ii}\big)^2-4\det[\f]}}2
\tag {6}\end{equation}

 Since we are assuming that $\det[\f]=0$, (6) becomes:
\begin{equation}q=\cases \sum_{i=1}^2f_{ii},\quad\mbox{or}\\0\endcases
\tag {7}\end{equation}

Choose $q=\sum_{i=1}^2f_{ii}$  and assume that
 $q\ne 0$, $q\neq f_{22}+f_{21}\beta_{12}$.
Then the second equation in the system implies that $b=0$ and
the system becomes:
$$\begin{array}{r@{\,}*{2}{c@{\:}}c@{\;=\;}c}
-f_{22}&a&+& f_{21}\,c &0\\
f_{12}&a&-&f_{11}\,c &0
\end{array}$$
 
If $f_{22}\neq 0$ then the above system has the general solution
$$ a=\frac {f_{21}}{f_{22}}\,c, \text{ where } c\in\C.$$ 
In particular, if we take $c=1$ then 
$p=\frac{f_{21}}{f_{22}}X_{12}+X_{22}$
satisfies
$$D\Big(\frac{f_{21}}{f_{22}}X_{12}+X_{22}\Big)=
\Big(\sum_{i=1}^2f_{ii}\Big)\Big(\frac{f_{21}}{f_{22}}X_{12}+X_{22}\Big).$$
On the other hand we also have that 
$$D(\det[\X])=\Big(\sum_{i=1}^2f_{ii}\Big)\det[\X].$$
Let
\begin{equation*}\theta =\frac{\frac{f_{21}}{f_{22}}X_{12}+X_{22}}{\det[\X]}.\end{equation*}
We have, 
 \begin{eqnarray*}\lefteqn{D(\theta)}\\&=&D\big(\frac{\frac{f_{21}}{f_{22}}X_{12}+X_{22}}
{\det[\X]}\big)\\
&=&\frac{D\big(\frac{f_{21}}{f_{22}}X_{12}+X_{22}\big)\det[\X]-
\big(\frac{f_{21}}{f_{22}}X_{12}+X_{22}\big)D(\det[\X])}
{\det[\X]^2}\\
&=&\frac{\big(\displaystyle\sum_{i=1}^2f_{ii}\big)
\big(\frac{f_{21}}{f_{22}}X_{12}+X_{22}\big)\det[\X]
-\big(\frac{f_{21}}{f_{22}}X_{12}+X_{22}\big)\big(\displaystyle\sum_{i=1}^2f_{ii}\big)\det[\X]}
{\det[\X]^2}\\
&=&0.\end{eqnarray*}
That is, $\theta$  \emph{is a new constant in} $F(\X)$.

Now we show that under the restrictions that we imposed on the
$\f$ it is possible to find a non-zero $f_{22}$. 

 Since we have a linear dependence relation among the $\X$, 
the  wronskian $W_1$  must be equal to zero. This
Wronskian can be expressed, up to a sign, as the following 
product of determinants:
$$ W_1=\left|\begin{array}{*{3}{c@{\quad}}c}
1&0&0&1\\f_{11}&f_{12}&f_{21}&f_{22}\\A&B
&E&F\\C&D&G&H\end{array}\right|
\left|\begin{array}{*{3}{c@{\quad}}c}
X_{11}&X_{12}&0&0\\X_{21}&X_{22}&0&0\\
0&0&X_{11}&X_{12}\\0&0&X_{21}&X_{22}\end{array}\right|
=M(\f)\det[\X]^2,$$
where
\begin{eqnarray*} A&=&f_{11}'+f_{11}^2+f_{12}f_{21},\\
B&=&f_{12}'+f_{11}f_{12}+f_{12}f_{22}\\
 C&=&f_{11}A+f_{21}B+A'\\
&=&3f_{11}f_{11}'+2f_{11}f_{12}f_{21}+2f_{12}'f_{21}+
f_{11}''+f_{12}f_{21}'+f_{11}^3,\\
D&=&f_{12}A+f_{22}B+B'\\
&=&2f_{11}'f_{12}+f_{11}^2f_{12}+f_{12}^2f_{21}+f_{21}f_{22}^2
+2f_{12}'f_{22}+f_{11}f_{12}'\\
&&\qquad +f_{12}f_{22}'+f_{12}''+f_{11}f_{12}f_{22},\\
E&=&f_{21}'+f_{21}f_{11}+f_{22}f_{21},\\
F&=&f_{22}'+f_{12}f_{21}+f_{22}^2,\\
G&=&f_{11}E+f_{21}F+E'\\
&=&2f'_{21}f_{11}+f_{21}f_{11}^2+f_{22}f_{21}f_{11}
+2f_{22}'f_{21}+f_{12}f_{21}^2\\
&&\qquad +f_{22}^2f_{21}+f_{21}''+f_{21}f_{11}'
+f_{22}f_{21}',\\ 
H&=&f_{22}F+f_{12}E+F'\\&=&f_{21}f_{11}f_{12}+2f_{22}f_{21}f_{12}+
+3f_{22}f'_{22}+2f_{12}f'_{21}+f'_{12}f_{21}\\
&&\qquad+f_{22}^3+f_{22}''.\end{eqnarray*} 
and $$M(\f)=\left|\begin{array}{*{3}{c@{\quad}}c}
1&0&0&1\\f_{11}&f_{12}&f_{21}&f_{22}\\A&B
&E&F\\C&D&G&H\end{array}\right|.$$

We have after simplifying using the hypothesis that $\det[\f]=0$,
\begin{equation*}\begin{split}M&(\f)=(f_{22}-f_{11})(f_{12}'f_{21}''-f_{21}'f_{12}'') 
+(f_{22}'-
f_{11}')(f_{12}''f_{21}-f_{12}f_{21}'') \\
&-f_{12}'f_{21}'(f_{11}-f_{22})^2-f_{12}f_{21}(f_{11}'-f_{22}')\\
 &+f_{12}f'_{21}(f_{11}f'_{11}+f_{22}f'_{22}-f_{11}'f_{22}-f_{11}f'_{22}
 +f_{22}''-f_{11}''
+f_{12}f_{21}'-f'_{12}f_{21})\\
&+f'_{12}f_{21}(f_{11}f'_{11}+f_{22}f'_{22}-f_{11}'f_{22}-f_{11}f'_{22}
 +f_{11}''-f_{22}'' 
 +f'_{12}f_{21}-f_{12}f'_{21}).\end{split}\end{equation*} 
Getting the above expression for $M(\f)$ took long and involved  computations.
We first computed the determinant directly and then we checked the result
using Dogson's method  \cite{D}, \cite{RR}.

The wronskian $W_1=0$ if and only if $M(\f)=0$. Now, observe that if
$f_{12}=0$ then $f_{12}'=0$ which implies that $B=0$ and $D=0$ as well.
Therefore $M(\f)=0$. So, if
we let $M(\Y)$ be the differential polynomial in the $\Y$
whose specialization to the $\f$ is $M(\f)$ then 
$M(\Y)$ is in the differential ideal 
\begin{equation*}\begin{split} \mathcal I&=\{\det[\Y],Y_{12}\}\\
&=\{Y_{11}Y_{22}-Y_{12}Y_{21},\,Y_{12}\}\\
&=\{Y_{11}Y_{22},\,Y_{12}\}\end{split}\end{equation*} 
of $C\{Y_{11},\,Y_{12},\,Y_{21},\,Y_{22}\}.$ It is easy to see that
$Y_{22}$ is not in $\mathcal I$. Indeed, suppose that 
\begin{equation}Y_{22}=p\,Y_{11}Y_{22}+q\,Y_{12}+r,\tag{8}\end{equation}
where $p,\,q\in C\{Y_{11},\,Y_{12},\,Y_{21},\,Y_{22}\}$,
$$r= \sum_{i,j}\big[p_i\,(Y_{11}Y_{22})^{(i)}+q_j\,Y_{12}^{(j)}\big]$$ with 
$p_i, q_j\in C\{Y_{11},\,Y_{12},\,Y_{21},\,Y_{22}\}$.
 
Now, consider the map
$$\psi :C\{Y_{11},\,Y_{21},\,Y_{22}\}\longrightarrow
C[Y_{11},\,Y_{21},\,Y_{22}]$$
given by  $\psi(Y_{22})=Y_{22}$ and 
$\psi(\Y)= 0$ for $i,\,j \ne 2$.
Let $\overline p=\psi(p)$, $\overline q=\psi(q)$, $\overline r=\psi(r)$.
We have that $\overline r=0$ and (8) becomes
$$Y_{22}=0. $$
which is impossible.$\hfill\qed$

\bibliographystyle{amsplain}

\begin{thebibliography}{10}

\bibitem{A-L} W.\   W. Adams and Ph. Loustaunau, \emph{An Introduction to Gr\"obner Bases},
Graduate Studies in Mathematics, American Mathematical Society (1994).


\bibitem{B-H} F.\    Beukers, G.\    Heckmann, Monodromy for the
hypergeometric function $_nF_{n-1}$, \emph{Invent. Math.} {\bf 95} (1989),
325--354.

 \bibitem{BB} A.\    Bialynicki-Birula, On the inverse problem of
Galois theory  of differential fields, \emph{Bull. Amer. Math. Soc.} {\bf
16} (1963), 960--964.


\bibitem{B} K.\    Boussel, \emph{Groupes de Galois des \'equations
hyperg\'eom\'etriques}, C. R. Acad. Sci. Paris. {\bf 309}, I (1989),
587--589.


\bibitem{D} C.\   L.\   Dodgson,
Condensation of determinants,
 \emph{Proc. Royal Soc. London} {\bf 15} (1866),  150--155.


\bibitem{dlr} A.\   Duval, M.\   Loday-Richaud, Kovacic's algorithm and
its applications to some families of special functions, \emph{AAECC Journal}
{\bf 3} (1992).

\bibitem{dm} A.\   Duval, C.\   Mitschi, Matrices de Stokes et groupe de
Galois des  \'equations hyperg\'eom\'etriques confluentes
g\'en\'eralis\'ees, \emph{Pacific J.  Math.} {\bf 138} (1989), 25--56.


\bibitem{goldman} L.\   Goldman, Specializations and the
Picard--Vessiot theory, \emph{Trans. Amer. Math. Soc.} {\bf 85} (1957),
327--356.

\bibitem{J} L.\   Juan,
 \emph{A Generic Picard--Vessiot extension with group $\GLn$ and the inverse differential 
Galois problem}, Ph.\  D.\   Thesis, University of Oklahoma  (2000).
 
 \bibitem{K} N.\   Katz, On the calculation of some differential Galois
groups, \emph{Invent. Math.} {\bf 87} (1987), 13--61.

\bibitem{K1}  J.\   Kovacic, The inverse problem in the Galois
theory of differential fields, \emph{Ann. of Math.} {\bf 89} (1969),
583--608.

\bibitem{K2}  J.\   Kovacic,  On the inverse problem in the Galois
theory of differential fields, \emph{Ann. of Math.} {\bf 93} (1971),
269--284.

\bibitem{M} A.\   Magid,
\emph{ Lectures in Differential Galois Theory}, 
University Lecture Series, American Mathematical Society  (1994).

\bibitem{miller} J.\    Miller  \emph{On Differentially Hilbertian
Differential Fields}, Ph.D. thesis, Columbia University (1970).

\bibitem{M1} C.\ Mitschi, Groupe de Galois des  \'equations
hyperg\'eom\'etriques confluentes g\'en\'eralis\'ees, \emph{C. R. Acad. Sci.
Paris.} {\bf 309}, I (1989), 217--220.

\bibitem{M2} C.\    Mitschi,  Differential Galois Groups of Generalized
Hypergeometric Equations: An Approach Using Stokes Multipliers,  \emph{Pacific J. Math.}
{\bf 176} 2, (1996), 365--405.


\bibitem{M-S2}  C.\    Mitschi and  M.\    F.\   Singer,
Connected Linear Groups as Differential 
Galois Groups,
\emph{Journal of  Algebra.} {\bf 184} (1996),
 333--361.

 
\bibitem{EN} E.\    Noether,
Gleichungen mit vorgeschriebener Gruppen, 
\emph{Math. Ann.} {\bf 78} (1918), 221--229.


\bibitem{P} M.\   van der Put,
Recent Work on Differential Galois Theory, 
\emph{S\'eminaire BOURBAKI}, 50\`eme ann\'ee,  (1997-1998), n$^o$ 849.

\bibitem{ramis1} J.\   -P.\    Ramis, About the solution of some inverse
problems in differential
Galois theory by Hamburger equations, in \emph{Differential Equations,
Dynamical Systems,
and Control Science}, Elworthy, Everett, and Lee, eds., Lecture Notes in
Pure and Applied
Mathematics. {\bf 157} (1994), 277--300.

\bibitem{ramis2} J.\   -P.\    Ramis,
 About the Inverse Problem in Differential Galois Theory: The Differential
Abhyankar Conjecture,
\emph{The Stokes Phenomenon and Hilbert's 16th Problem};
 B.\   I.\   J.\     Braaksma, et al., eds.\  World Scientific,
 Singapore (1996)
 
\bibitem{RR} D.\   P.\    Robins and H. Rumsey, Jr.,
Determinants and Alternating Sign Matrices,  
\emph{Advances in Mathematics.} {\bf 62} (1986),  169--184.

\bibitem{S}  M.\   F.\   Singer,
 Direct and Inverse Problems in Differential Galois Theory,
\emph{ Selected Works of Ellis Kolchin with Commentary},
 Bass, Buium, Cassidy eds., American Mathematical Society (1999),
527--554.

\bibitem{S2} M.\   F.\   Singer, Moduli of linear differential equations
on the Riemann
sphere with fixed Galois groups, \emph{Pacific J.  Math.} {\bf 106}, 2 (1993),
343--395.

\bibitem{su} M.\   F.\   Singer, F.\   Ulmer,  Liouvillian and algebraic
solutions of second and third order linear differential equations, \emph{Journal
of Symbolic Computation.}  {\bf 16} (1993), pp. 37--74.


\bibitem{TT} C.\   Tretkoff and M.\   Tretkoff, {\em  Solution of the inverse
problem of
differential Galois theory in the classical case}, \emph{Amer. J. Math.} {\bf
101} (1979), 1327--1332.

\bibitem{uw} F. Ulmer, J.\   -A.\   Weil, Note on Kovacic's algorithm,
\emph{J. Symbolic Comput.} {\bf 28}, 2, (1996) 179--200.

\bibitem{W}  J.\   -A.\    Weil, \emph{Constantes et polyn\^omes de Darboux en alg\`ebre diff\'erentielle:
applications aux syst\`emes diff\'erentiels lin\'eaires}, 
Ph.D. Thesis,  \'Ecole Polytechnique de France (1995).
\end{thebibliography}

\end{document}